%% file: FINAL_VERSION.tex
\documentclass[journal]{IEEEtran}    
\input{preamble.tex}

\title{Buildings-to-Grid Integration Framework}
\author{Ahmad F. Taha,~\IEEEmembership{Member,~IEEE, }Nikolaos Gatsis,~\IEEEmembership{Member,~IEEE,} Bing Dong,~\IEEEmembership{Member,~IEEE,}\\ Ankur Pipri,~\IEEEmembership{Student Member,~IEEE,} and Zhaoxuan Li,~\IEEEmembership{Student Member,~IEEE.}
	\thanks{
Manuscript received March 18, 2017; revised June 18, 2017 and August 26, 2017; accepted October 1, 2017. This work was supported by the National Science Foundation under Grants CBET-1637249 and ECCS-1462404. Paper no. TSG-00379-2017.

Ahmad F. Taha, Nikolaos Gatsis, and Ankur Pipri are affiliated with the Department of Electrical and Computer Engineering, and Bing Dong and Zhaoxuan Li are affiliated with the Department of Mechanical Engineering, all at
	The University of Texas at San Antonio. Emails: 
			\{ahmad.taha,nikolaos.gatsis,ankur.pipri,bing.dong,zhaoxuan.li\}@utsa.edu. }}
\begin{document}
\maketitle
\begin{abstract}
This paper puts forth a mathematical framework for Buildings-to-Grid (BtG) integration in smart cities. The framework \textit{explicitly} couples power grid and building's control actions and operational decisions, and can be utilized by buildings and power grids operators to simultaneously optimize their performance. Simplified dynamics of building clusters and building-integrated power networks with algebraic equations are presented---both operating at different time-scales.  A model predictive control (MPC)-based algorithm that formulates the BtG integration and accounts for the time-scale discrepancy is developed. The formulation captures dynamic and algebraic power flow constraints of power networks and is shown to be numerically advantageous. The paper analytically establishes that the BtG integration yields a reduced total system cost in comparison with decoupled designs where grid and building operators determine their controls separately. 
The developed framework is tested on standard power networks that include thousands of buildings modeled using industrial data.   Case studies demonstrate building energy savings and significant frequency regulation, while these findings carry over in network simulations with nonlinear power flows and mismatch in building model parameters. Finally, simulations indicate that the performance does not significantly worsen when there is uncertainty in the forecasted weather and base load conditions.
\end{abstract}

\section{Introduction, Prior Art, Paper Contributions}~\label{sec:Intro}

\IEEEPARstart{B}{y} 2050, a staggering 70\% of the world's population is bound to live and work in cities \cite{Pop2050}. A recent assessment from the World Bank suggested that two-thirds of global energy consumption can be attributed to cities, leading to 71\% of global direct energy-related greenhouse gas emissions~\cite{swilling2013city}. Smart cities consist of sustainable and resilient infrastructures, where buildings are a major constituent. Building energy consumption contributes to more than 70\% of electricity usage---profoundly impacting power grid's operation. Futuristic cities equipped with optimized building designs have the auspicious potential to play a pivotal role in reducing global energy consumption while maintaining stable electric-grid operations. As buildings are physically connected to the electric power grid, it is natural to understand their coupling and develop a framework for Buildings-to-Grid (BtG) integration. To understand the role and impact of BtG integration, the authors in~\cite{LAWRENCE2016273} provide relevant research questions for BtG integration. 

The installation of smart meters in buildings and across the power grids  enables the BtG integration, which can transform passive buildings into active dispatchable demand resources. The U.S.\ Department of Energy has highlighted the multiple benefits and opportunities of BtG integration~\cite{DoE-BtGreport}: 1) Buildings can enjoy significant energy savings; 2) the grid's resources are more efficiently utilized, as peak demand is curbed; 3) the grid can become more stable with fewer frequency excursions; 4) the need for bulk generation and transmission investments is deferred; and 5) with a BtG integration platform in place, distributed energy resources at the buildings' premises, such as photovoltaic units and electric vehicles, can be more efficiently integrated with the power grid, turning into significant assets. 
Besides the aforementioned technical motivating factors, the mathematical intuition behind the BtG integration is that when the grid and the buildings jointly optimize their control decisions, they have the potential to yield larger system benefits than when they make these decisions separately. This paper aspires to develop a BtG integration framework, and make this mathematical intuition precise.

Various studies address a breadth of computational and experimental aspects of BtG integration. An overview of demand response potential from smart buildings is presented in~\cite{qi2016demand}. An experimental architecture that enables smart buildings is proposed in~\cite{stamatescu2016building} with a focus on heating, ventilation, and air conditioning (HVAC) systems and grid integration.  A bi-level optimization framework for commercial buildings integrated with a distribution grid is proposed in~\cite{Razmara2016}.  Detailed dynamical models for buildings with multiple zones (upper level) and an operational model for the distribution grid with voltage and current balance equations (lower level) are included; nonetheless, a dynamical model of the power grid, suitable for modeling frequency excursions, is missing. 

The regulation service provision by smart buildings is investigated in~\cite{Bilgin2016}, where price signals are exchanged between grid and building operators to alter building energy consumption. Other BtG integration studies have shown that grid-aware building HVAC controls  can provide frequency regulation or other ancillary services to the grid \cite{zhao2013evaluation, Zhao20151,Blum20141,hao2014ancillary,Lin2017}, largely without sacrificing the occupants' comfort. The load-shifting capability of buildings has also been explored~\cite{7315254}. Explicit account of the grid dynamics and power flows is on the other hand missing from the previously mentioned works.

HVAC controls and building dynamics are typically modeled as linearized dynamical systems~\cite{Wang20061927}. The modeling particulars depend on the size of the building cluster. As the number of buildings involved in the analysis increase, the dynamical models tend to become simpler---for obvious computational purposes.
A typical thermal resistance and capacitance circuit model can be used to represent heat transfer and thermodynamical properties of the building envelope, and is widely used in building control studies \cite{kummert2011using,yahiaoui2006model,dong2010integrated,maasoumy2014modeling}. Given these models for building dynamics, various control routines have been developed for building controls.
Currently, many commercial buildings use PID controllers for HVAC systems~\cite{Afram2014343}. However, model predictive control (MPC) has proved to be advantageous with respect to PID controllers~\cite{morocsan2011distributed,patel2016distributed,Mirakhorli2016499}.  

Building HVAC control via MPC has been investigated under a wide range of scenarios and setups. For example, uncertainty in the building MPC formulation is considered in~\cite{Oldewurtel2010,rostampour2016probabilistic,oldewurtel2012use}. Centralized building MPC routines are proposed in \cite{Ma2012,Chen2013}, and the works in~\cite{morocsan2011distributed,patel2016distributed,standardi2015economic} investigate decentralized or distributed solvers to building MPC problems.  In addition, explicit MPC routines have also been developed in the context of building control studies~\cite{koehler2013building}. Other works focus on integrating occupancy behavior and its impact on indoor temperature variations, while still attempting to obtain optimal control laws~\cite{Mirakhorli2016499,dobbs2014model}. The majority of the aforementioned works show significant energy savings given different system dynamics, forecast and parametric uncertainty, and computational limitations. 

While the aforementioned research investigated different challenging problems related to BtG integration and building MPC routines, none of these studies produces a high-level mathematical framework that buildings and power grids operators can simultaneously utilize to optimize their performance---a framework that \textit{explicitly} couples power grid and building control actions and operational decisions. In addition, the majority of the previously mentioned studies focus on one or a group of buildings and the corresponding impact on the power grid, rather than clusters of thousands of buildings in smart cities. The main challenges associated with creating such a framework that addresses the aforementioned research gaps are as follows.
\begin{itemize}[leftmargin=*]
	\item Building control systems are neither connected to each other, nor integrated with the grid. Consequently, a unified optimal energy control strategy---even if it is decentralized---cannot be achieved unless there is a framework that facilitates this integration, in addition to the willingness of building operators to contribute to this framework.
	\item Grid and building dynamics and control actions clearly operate at two different time-scales. While the grid controls and states are often in seconds,  the building state dynamics and controls are much slower, often in minutes. Coupling the two dynamic systems together entails addressing this time-scale discrepancy.
	\item Existence of algebraic equations in grid dynamics, resulting in differential algebraic equations (DAE), coupled with the different time-scales, complicates modeling and analysis of BtG integration. In fact, these algebraic equations depict the interdependence between grid and building dynamics.
\end{itemize}
The chief contribution of this paper is a novel mathematical framework for BtG integration that addresses the aforementioned challenges in a structured and principled way. The MPC-based framework couples building dynamics, grid dynamics that include the network frequency, and the power flow equations. The objective  is to generate local control actions for buildings and power generators such that the overall performance is optimized in terms of stability, energy savings, and other socio-economic metrics. The time scale discrepancy between the grid and building dynamics is explicitly accounted for in the developed MPC-based optimization formulation. The building-integrated power network dynamics are modeled by DAEs. In order to include the DAEs into the optimization framework, the DAEs are discretized using Gear's method. While this discretization method has been the basis for the numerical solution of DAEs for several decades~\cite{sincovec1981analysis}, it is the first time that it is brought in to facilitate the development of a BtG integration framework. The paper analytically establishes that the BtG integration yields a reduced total system cost with respect to decoupled designs in which grid and building operators determine their controls separately. 

 The developed framework is tested in standard IEEE networks that include hundreds to thousands of buildings modeled using ASHRAE data. The building HVAC load is driven by realistic ambient weather patterns and typical temperature requirements for commercial buildings. The simulations are performed for various power networks, and a reduction of up to 20\% and 43\% in total system cost---with respect to two other decoupled designs---is demonstrated. Finally, simulations indicate that the performance does not significantly worsen when there is uncertainty in the building parameters or forecasted weather and base load conditions.

The remainder of this paper is organized as follows. In Sections~\ref{sec:buildingsSection} and~\ref{sec:PowerNetwork}, we present the dynamics of building clusters and of the building-integrated power network. Optimal power flow is also integrated in these models. In Section~\ref{sec:OCP}, we propose our approach for BtG integration, while addressing the aforesaid challenges. Then, the optimal control problem that models the BtG integration is formulated. A customized algorithm is also developed to seamlessly include the optimal power flow into the integrated framework. Section~\ref{sec:anlytical} produces an analytical discussion on the advantages of the BtG framework over decoupling the optimization of buildings and power grids. Case studies with realistic building parameters and grid constraints are given in Section~\ref{sec:results}. Finally, future work is outlined in Section~\ref{sec:summary}.

\section{Building Clusters Dynamics}~\label{sec:buildingsSection}
The patterns of energy usage in buildings are impacted by local climate, heat transfer through the building envelope, daily operation, and occupancy behaviors. Detailed energy models have been developed based on physics and statistics to simulate heat transfer in buildings. 

For a large-scale application such as BtG integration, it is unrealistic to consider every thermal zone of each building---thousands of buildings will generate millions of zones. This approach would produce a highly intractable BtG integration problem.  Hence, at a BtG integration level, the amount of cooling energy needed which is optimized to minimize the total operation cost is allocated to each building; we define this quantity as $P_{\mathrm{HVAC}}\itn{l}$ for building $l$. Then, at the local level of building $l$, the decision variables of the air-side system (setpoints for air-handling units, damper opening for terminal systems) and the water-side system (flows for pumps, chiller temperature setpoints) can be optimized to maintain the preferred zone temperature, while not exceeding the cooling load limits set by $P_{\mathrm{HVAC}}\itn{l}$. This approach of solving for the cooling loads and then feeding the setpoints to local lower level problems is common in recent building studies; see~\cite{patel2016distributed}. 

\begin{figure}[t]
	\centering
	\includegraphics[scale=0.225]{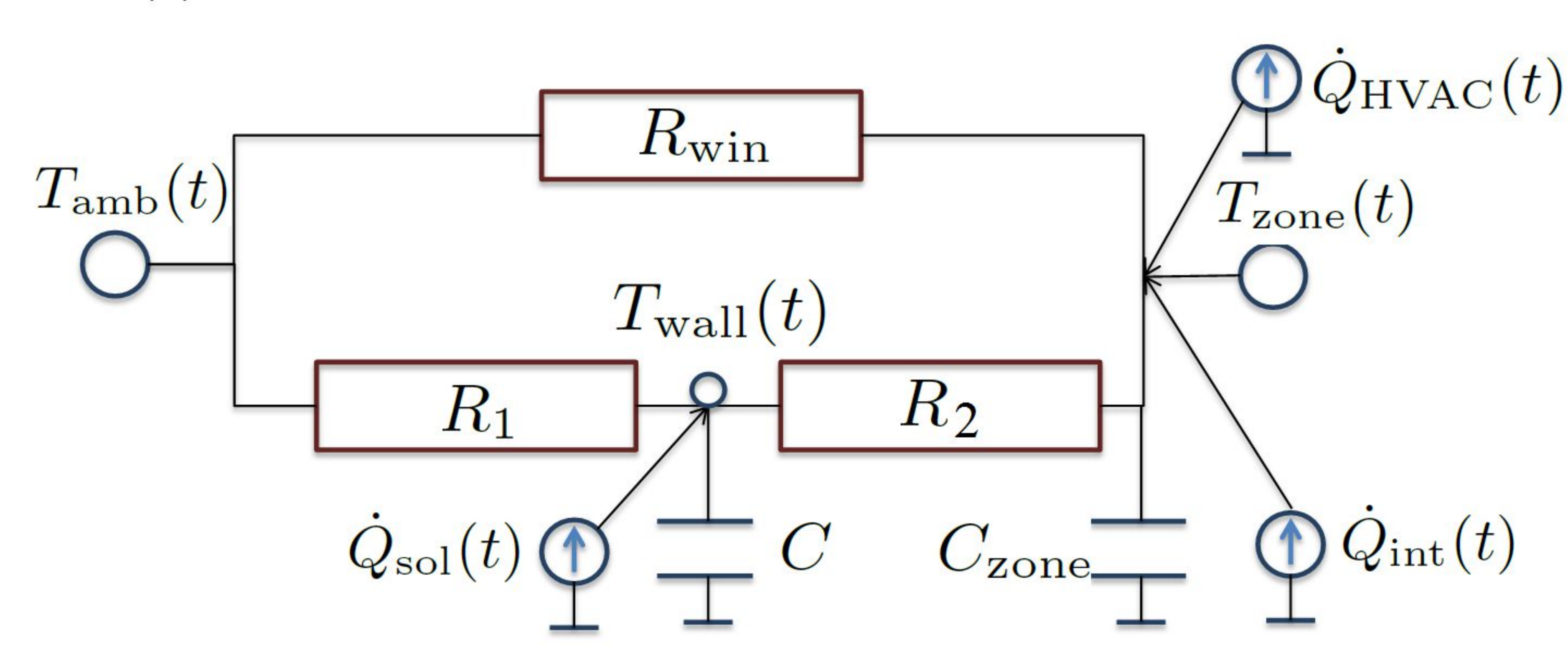}
	\caption{An RC-network model for a thermal zone~\cite{Dong2013}.}
	\label{fig:buildingmodel} 
\end{figure}

In this paper, we focus on the high-level commercial buildings problem described above. We use a typical thermal resistance and capacitance (RC) network to model heat transfer and the thermodynamics of the building envelope, which has been widely used in building control studies~\cite{kummert2011using,yahiaoui2006model,dong2010integrated,maasoumy2014modeling}.  The RC network model assumes a steady-state heat transfer through the building envelope. Considering that the building dynamics have time constant of hours, this model is sufficient for a high-level BtG integration study. A typical three-resistance and two-capacitance (3R-2C) model is shown in Fig.~\ref{fig:buildingmodel}. 
In this model, a building is treated as a \textit{super-zone} where the resistance parameters represent the thermal resistance of the building structure, the external facades' convection, and the internal walls' convection.  Building dynamics with temperatures $T_{\mathrm{wall}}(t)=T_{\mathrm{wall}}$ and $T_{\mathrm{zone}}(t)=T_{\mathrm{zone}}$ are written as
 \begin{eqnarray*} 
\dot{T}_{\mathrm{wall}}  & =&  \dfrac{T_{\mathrm{amb}}  -T_{\mathrm{wall}}  }{CR_2}+\dfrac{T_{\mathrm{zone}}  -T_{\mathrm{wall}} }{CR_1}+\frac{\dot{Q}_{\mathrm{sol}} }{C}  \label{equ:buildingDynamics} \\
\dot{T}_{\mathrm{zone}}  & = & \dfrac{T_{\mathrm{wall}}  -T_{\mathrm{zone}}  }{C_{\mathrm{zone}}R_1}+\dfrac{T_{\mathrm{amb}}  -T_{\mathrm{zone}} }{C_{\mathrm{zone}}R_{\mathrm{win}}} +\frac{\dot{Q}_{\mathrm{int}}  +  \dot{Q}_{\mathrm{HVAC}} }{C_{\mathrm{zone}}}, \nonumber 
\end{eqnarray*}
where $R_{\mathrm{win}}$, $R_2$, and $R_1$ are physical parameters of the building envelope;
 $C$ is the lumped thermal capacity of all walls and the roof;
$C_{\text{zone}}$ is the thermal capacity of the zone;
 $\dot{Q}_{\mathrm{sol}}(t)$ is the total absorbed solar radiation on the external wall;
$\dot{Q}_{\mathrm{int}}(t)$ is the total internal heat gain from space heat sources such as desktops, people, and lights; $T_{\mathrm{zone}}(t)$ and $T_{\mathrm{wall}}(t)$ are respectively the zone (space) and wall temperatures; and
$T_{\mathrm{amb}}(t)$ is the outside ambient temperature. 
The cooling load can be calculated as
		$\dot{Q}_{\mathrm{HVAC}}(t)=\mu_{\mathrm{HVAC}}P_{\mathrm{HVAC}}(t),$ 
		where	$P_{\mathrm{HVAC}}$ is the power consumed by the HVAC system, and
		 $\mu_{\mathrm{HVAC}}$ is the coefficient of performance of the HVAC system.
The dynamics of building $l$ are described by
\begin{equation}~\label{equ:BuildingSS}
	\dot{\m x}_b\itn{l}(t) = \m A_b\itn{l} \m x_b\itn{l}(t) + \m B_{u_{b}}\itn{l} \m u_b\itn{l}(t) + \m B_{w_{b}}\itn{l} \m w_b\itn{l}(t),
\end{equation}
where $\m x_b\itn{l}  = [T_{\mathrm{wall}} \; T_{\mathrm{zone}}]^{\top}_l$ is the state of building $l$; $\m u_b\itn{l} = \bmat{P_{\mathrm{HVAC}}}_l$ is the control input variable; $\m 	w_b\itn{l} = [ T_{\mathrm{amb}} \;  \dot{Q}_{\mathrm{sol}}  \; \dot{Q}_{\mathrm{int}}]^{\top}_l$ is a random uncontrollable input. Various methods have been developed to provide an estimate of $\m w_b\itn{l}$, denoted by $\hat{\m w}_b\itn{l}$, for each building~\cite{dong2014real}. The system state-space matrices in~\eqref{equ:BuildingSS} are defined as
\begin{eqnarray*}
	\m A_b\itn{l}  &=&  \bmat{
		-\dfrac{1}{C}\left(\dfrac{1}{R_1}+\dfrac{1}{R_2}\right) & 		\dfrac{1}{CR_1}\\
		\dfrac{1}{C_{\mathrm{zone}}R_1} & 	-\dfrac{1}{	C_{\mathrm{zone}} }\left(\dfrac{1}{R_1}+\dfrac{1}{R_{\mathrm{win}}}\right) 
	}_l \\
	\m B_{u_{b}}\itn{l}  &=&   \bmat{
		0\\
		\dfrac{\mu}{C_{\mathrm{zone}}}}_l ,
	\m B_{w_{b}}\itn{l} =   \bmat{
		\dfrac{1}{CR_2}	& \dfrac{1}{C} & 0 \\
		\dfrac{1}{C_{\mathrm{zone}}R_{{\mathrm{win}}}} 	& 0 &  \dfrac{1}{C_{\mathrm{zone}}}}_l.
\end{eqnarray*}
The notation $[\cdot]_l$ implies that each building $l$ has a different set of RC parameters. In this paper, we consider clusters of buildings with each cluster connected to a power grid node. Since we aim to understand the impact of buildings' contribution to frequency regulation and overall energy consumption costs, we present the dynamics of building clusters
\begin{equation}~\label{equ:CBuildingDynamincs2}
\dot{\m x}_b(t) = \m A_b \m x_b(t) + \m B_{ u_{b}} \m u_b(t) + \m B_{w_{b}} \m w_b(t)\,,\,\,
	\end{equation}
where $n_b$ is the total number of buildings in the network;
 $\m x_b \in \mathbb{R}^{2n_b}, \m u_b \in \mathbb{R}^{n_b},$ and  $\m w_b \in \mathbb{R}^{3n_b}$. In the absence of communication between buildings, the state-space matrices $\m A_b$, $\m B_{u_{b}}$, and $\m B_{w_{b}}$ will all be block diagonal matrices:
$\m A_b = \diag(\m A_b\itn{1},\ldots, \m A_b\itn{n_{b}}) \in \mathbb{R}^{2n_b \times 2n_b},\; \m B_{u_{b}}=\diag(\m B_{u_{b}}\itn{1},\ldots, \m B_{u_{b}}\itn{n_b}) \in \mathbb{R}^{2 n_b \times n_b}, \; \m B_{w_{b}} = \diag(\m B_{w_{b}}\itn{1},\ldots, \m B_{w_{b}}\itn{n_b}) \in \mathbb{R}^{2 n_b \times  3n_b}.$
\begin{rem}\normalfont 
Here, we assume that the variables to be solved for are the high-level, total cooling load setpoints for individual buildings $\m u_b(t)$. From this cooling load, building operators solve low-level control and optimization problems. This is customary in building control studies; see for comparison~\cite{patel2016distributed}.
\end{rem}
\section{Building-Integrated Power Network: Dynamics And Connection to OPF}\label{sec:PowerNetwork}
In this section, we present the dynamics of the building-integrated power network and define the main variables involved in the BtG integration framework.  In addition, we discuss the connection of the BtG integration model with the optimal power flow (OPF).
\subsection{DAE Dynamics of a Power Network with Building Loads}~\label{sec:LTIDynamics}
 Let $\mathcal{B}=\{1,\ldots,n\}$ and $\mathcal{G}=\{1,\ldots,n_g\}$ denote the sets of buses and generators in a power network. Also, let 
$\mathcal{N}_k$ be the neighborhood set of adjacent nodes connected to the $k$th bus.  Generators are indexed by $m\in\mathcal{G}$. The mechanical input power to the $m$th generator is denoted by $P_m$. Define generator-to-node and building-to-node incidence matrices $\m \Gamma \in \mathbb{R}^{n \times n_g}$ and $\m \Pi \in \mathbb{R}^{n \times n_b}$ with entries given by
\begin{subequations}
\label{eq:inddef}
\begin{align}
\gamma_{k,m} &= \left\{ 
	\begin{array}{l l l} 
		 1 & \hspace{0.2cm}\text{if} & \text{generator $m$ is attached to bus}\,\,k \\
		 0 &\hspace{0.2cm} \text{} & \text{otherwise},
	\end{array} \right.   \label{equ:gammadef}\\
	\pi_{k,l} & = \left\{ 
		\begin{array}{l l l} 
			1 & \hspace{0.2cm}\text{if} & \text{building $l$ is attached to bus}\,\, k \\
			0 &\hspace{0.2cm} \text{} & \text{otherwise}.
		\end{array} \right.   \label{equ:pidef}
\end{align}
\end{subequations}
The transients of the $k$th bus in a power network can be modeled by the swing equation which relates the rotor angle $\delta$
with the angular velocity $\dot{\delta}$ and the angular acceleration $\ddot{\delta}$~\cite{Taylor2015}. Define $M_k$ and $D_k$ as the inertia and damping coefficients of the generator located $k$th bus; if the $k$th bus does not have a generator, then $M_k=0$ and $D_k=0$.
The swing equation for the $k$th bus ($k \in \mathcal{B}$) can be written as
\begin{eqnarray}
	M_k\ddot{\delta}_k(t)+D_k\dot{\delta}_k(t)&=&\gamma_{k,m}P_{m}(t)-P_{L_k}(t) \nonumber \\
	& & -\displaystyle \sum_{j \in \mathcal{N}_k} b_{kj} \sin(\delta_k(t)-\delta_j (t)).~\label{equ:SwingEquation2}
\end{eqnarray}
{The load at bus $k$, $P_{L_k}(t)$ is described by
\begin{equation}
P_{L_k}(t)= P_{\mathrm{BL}_k}(t) +  D_k'\dot{\delta}_k(t)  + \sum_{l=1}^{n_b} \pi_{k,l}P_{\mathrm{bldg}}\itn{l}(t).
\label{equ:busload}
\end{equation} 
In~\eqref{equ:busload}, the first two terms represent uncontrollable loads, while the last one defines the controllable load. Specifically, $P_{\mathrm{BL}_k}(t)$ denotes the frequency-insensitive uncontrollable base load at bus $k$, which is typically available via forecasts. The term  $D_k'\dot{\delta}_k$ denotes the frequency-sensitive portion of the uncontrollable load at bus $k$. It is assumed specifically that a portion of the load at bus $k$ responds linearly to the frequency variations, which is a classical model~\cite{OSullivan96}; the linear coefficient is $D_k'$, and the frequency is the derivative of the angle $\dot{\delta}_k(t)$. It must be emphasized that the uncontrollable loads are not optimization variables. The term  $\sum_{l}\pi_{k,l}P\itn{l}_{\mathrm{bldg}} (t)$ defines the load from buildings indexed by $l$ and attached to bus $k$ participating in regulation. The building load is further decomposed as $P\itn{l}_{\mathrm{bldg}} (t)=P_{\mathrm{HVAC}}\itn{l}(t)+P_{\mathrm{misc}}\itn{l}(t),$
where $l \in \{1,2,\ldots, n_b\}$ is the index of buildings.
 The quantity $P_{\mathrm{HVAC}}\itn{l}(t)$ denotes the portion of controllable power consumption of building $l$, while 
$P_{\mathrm{misc}}\itn{l}(t)$ represents the uncontrollable miscellaneous power consumption of building $l$ such as lighting, computers, equipment, elevators---amounting to a building's base-load. The quantity $P_{\mathrm{HVAC}}\itn{l}(t)$ is an optimization variable, and $P_{\mathrm{misc}}\itn{l}(t)$ is typically available via forecasts.}

We can now rewrite \eqref{equ:SwingEquation2} as
\begin{eqnarray}
 M_k\ddot{\delta}_k(t) &=& - D_k\dot{\delta}_k (t) + \gamma_{k,m}P_{m}(t) \nonumber \\
 && \hspace{-0.4cm} -\sum_{j \in \mathcal{N}_k} b_{kj} \sin\left(\delta_k (t)-\delta_j(t) \right) -D_k'\dot{\delta}_k (t) \nonumber \\
&& \hspace{-0.4cm}-P_{\mathrm{BL}_k}(t)  - \sum_{l=1}^{n_b}\pi_{k,l}\left( P_{\mathrm{HVAC}}\itn{l}(t) +P_{\mathrm{misc}}\itn{l}(t) \right).~\label{equ:SwingEquation3}
\end{eqnarray}
In~\eqref{equ:SwingEquation3}, $P_m(t)$ for generator bus $m$ can be written as
$ P_m(t)=\bar{P}_m+\Delta P_m(t),$
where $\bar{P}_m$ is a solution of an optimal power flow problem---computed every 15 minutes---and $\Delta P_m(t)$ is the deviation from the setpoint $\bar{P}_m$, which will be furnished by the proposed BtG integration framework.

The angular velocity is $\dot{\delta}_k=\omega_k$, where $\omega_k=\omega^{\mathrm{true}}_{k}-\omega_0$, $\omega^{\mathrm{true}}_k$ is the actual frequency of the $k$th bus, and $\omega_0$ is the synchronous frequency, e.g., $2\pi 60 \,\,\mathrm{rad/sec}$ in North America.
Given~\eqref{equ:SwingEquation3}, we obtain two first-order differential equations representing the dynamics of the $k$th bus
\begin{eqnarray*}
	\dot{\delta}_k(t) &=&\omega_k(t)\label{equ:sw1}\\
	M_k\dot{\omega}_k(t)&=&\displaystyle - \left({D_k+D_k'}\right)\omega_k(t) +\gamma_{k,m}{P_{m}(t)}\label{equ:sw2}\\
& & - {P_{B_k}(t)}  - \sum_{l=1}^{n_b}\pi_{k,l}\left( P_{\mathrm{HVAC}}\itn{l}(t)+P_{\mathrm{misc}}\itn{l}(t)\right)\nonumber  \\
	& &-\sum_{j\in \mathcal{N}_k} {b_{kj}} \sin \left( \delta_k(t)-\delta_j(t)\right),\, \, k \in \mathcal{B}.\nonumber
\end{eqnarray*}
The resulting state-space model is a nonlinear system, and is formulated as
\begin{eqnarray}
 \m E_g\dot{\m x}_g(t) & =& \m A_g\m x_g(t)+ \pmb \Phi(\pmb \delta (t))+\m A_{u_{b}} \m u_{b} (t)\nonumber \\
 &&+\m B_{u_{g}} \m u_{g}(t)+\m B_{w_{g}}\m w_{g}(t),\label{equ:PlantSS} 
\end{eqnarray}
where
	 $\m x_g(t)  = [\delta_1 \; \ldots \; \delta_n \;  \omega_1 \;  \ldots \; \omega_n]^{\top}=[\pmb{\delta}^{\top}(t) \; \pmb{\omega}^{\top}(t)]^{\top}$ is the state of the \textbf{g}rid; $\pmb \Phi(\pmb \delta (t))$ is the vectorized nonlinear power flow equations in~\eqref{equ:SwingEquation3};
	 $\m u_b (t) = [P_{\mathrm{HVAC}}\itn{1} \; \ldots \; P_{\mathrm{HVAC}}\itn{n_b}]^{\top}$ is the control input vector of the buildings, as defined in~\eqref{equ:CBuildingDynamincs2}, and $\m u_{g}(t) =\bar{\m u}_g+ \Delta \m u_g(t)
=[\bar{P}_1+\Delta P_1(t) \;\ldots \;\bar{P}_{n_{g}}+\Delta P_{n_{g}}(t)]^{\top}$ is the power network's control variable;
	$\m w_{g}(t) =[\m w_{\mathrm{BL}}^{\top}, \m w_{\mathrm{misc}}^{\top}]^{\top} = [P_{\mathrm{BL}_{1}}\;\ldots\;P_{\mathrm{BL}_{n}}, \, P_{\mathrm{misc}}\itn{1} \; \ldots \; P_{\mathrm{misc}}\itn{n_b} ]^{\top}$
	is a random vector collecting the nodal base loads and miscellaneous building loads. Load forecasting is a very mature area; in the sequel, the forecast of $\m w_{g}(t)$, denoted by $\hat{\m w}_{g}(t)=[\hat{\m w}_{\mathrm{BL}}^{\top}, \hat{\m w}_{\mathrm{misc}}^{\top}]^{\top}$, is assumed to be available.
	The state-space matrices in~\eqref{equ:PlantSS} are obtained as follows
\begin{align*}
	 \m E_g & =  \bmat{ 
	 \m I_n & \m 0_{n\times n}  \\
	 \m	0_{n\times n}   &   \m M },
	\m A_g  =   \bmat{ 
	\m 0_{n\times n}  &  \m I_n \\
			\m 0_{n\times n}&  -  \m D\\
	}, \m A_{u_{b}}  =  \bmat{	\m 0_{n\times n_b} \\  -\m\Pi}\\
	\m M &=\diag(M_1,\ldots,M_n), \, \bm{\Phi}(\bm{\delta})= [ \bm{0}_n^\top, \Phi_1(\bm{\delta}), \dots, 	\Phi_n(\bm{\delta})]^\top \\
	\Phi_k& =  \sum_{j\in \mathcal{N}_k} {b_{kj}} \sin \left( \delta_k-\delta_j\right), \; k=1,\ldots,n \\
	\m D&=\diag(D_1+D_1',\ldots,D_n+D_n'),\;  \m B_{u_{g}}   =   \bmat{	\m 0_{n\times n_b}  \\  \m \Gamma \\} \\
	\m B_{w_{g}} &=\diag(\m B_{\mathrm{BL}},\m B_{\mathrm{misc}}),	
	 	 \m B_{\mathrm{BL}}  =  \bmat{	\m 0_{n\times n} \\  -\m I_n \\},
	 	 	 	 \m B_{\mathrm{misc}}  =   \bmat{	\m 0_{n\times n_b} \\  -\m \Pi}
\end{align*}
\noindent where 
$\m E_g \in \mathbb{R}^{2n\times 2n}$ is a singular matrix,
$\m A_g \in \mathbb{R}^{2n\times 2n}$, $\m A_{u_{b}}  \in \mathbb{R}^{2n\times n_b}$,
$\m B_{u_{g}} \in \mathbb{R}^{2n\times n_g}$, and $\m B_{w_{g}} \in \mathbb{R}^{2n\times (2n+n_b)}$. 

Dynamic systems of the form~\eqref{equ:PlantSS} are called differential algebraic equations (DAE), and the systems following such equations are called descriptor systems. The chief difference between descriptor systems and standard dynamical systems is that the former  have a matrix $\m E_g$ that multiples $\dot{\m x}_g$ and may be singular. This is indeed the case here, where entries in the $\m M$ matrix corresponding to non-generator buses are zero, giving rise to entire rows of zeros. The next section describes how the generation setpoints in~\eqref{equ:PlantSS} are computed. 
\subsection{Connection to the Optimal Power Flow}
Recall that $\m u_g(t)$ is written as  $\m u_g(t)=\bar{\m u}_g+\Delta \m u_g(t),$
where $\bar{\m u}_g$ contains the generator setpoints, and $\Delta \m u_g(t)$ is the real-time deviation from these setpoints that automatically drives the power grid to stability after load deviations. Typically, the setpoints are computed every 5--15 minutes through solving economic dispatch or OPF routines~\cite{Taylor2015}. A linearized OPF (LOPF) problem is described by the following program
	\begin{subequations}\label{equ:DCOPF}
\begin{eqnarray}
\hspace{-0.3cm}\textbf{LOPF:}\;\;\;\underset{\bar{\m u}_g=\{\bar{u}_{g_{i}}\}_{i=1}^{n_g}}{\minimize} & & \hspace{-0.2cm}J(\bar{\m u}_g)=\bar{\m u}_g^{\top}\m J_{u_g} \bar{\m u}_g + \m b_{u_{g}}^{\top}\bar{\m u}_g+c_{u_{g}}~\label{equ:DCOPF0} \\ 
 \mathrm{subject\;to} & & \bar{\m u}_g^{\mathrm{min}} \leq \bar{\m u}_g \leq \bar{\m u}_g^{\mathrm{max}}   ~\label{equ:DCOPF1} \\
       & & \hspace{-2.4cm}\left(\m \Gamma\bar{\m u}_g-\m \Pi (\m u_{{b}}+\hat{\m w}_{\mathrm{misc}})- \hat{\m w}_{\mathrm{BL}} \right)^{\top} \m 1_n =0~\label{equ:DCOPF2} \\
         & &\hspace{-3cm} |\m L_{\mathrm{ptdf}} \left(\m \Gamma\bar{\m u}_g-\m \Pi (\m u_{{b}}+\hat{\m w}_{\mathrm{misc}})- \hat{\m w}_{\mathrm{BL}} \right) | \leq \m F^{\mathrm{max}}  ~\label{equ:DCOPF3}.
\end{eqnarray}
	\end{subequations}
{In~\eqref{equ:DCOPF}, $J(\bar{\m u}_g)$ is a convex cost function that represents the generators' cost curves. Constraint~\eqref{equ:DCOPF1} represents the safety upper and lower bounds on the generator's active power.  Vectors  $\m u_b$,  ${\m w}_{\mathrm{misc}}$, and ${\m w}_{\mathrm{BL}}$ were introduced in the previous subsection and represent respectively the building HVAC loads, building miscellaneous loads, and nodal base loads. We use the notations $\hat{\m w}_{\mathrm{misc}}$ and $\hat{\m w}_{\mathrm{BL}}$  to emphasize that the respective forecasted versions of the building miscellaneous loads and nodal base loads enter the OPF.  
	
Matrix $\bm\Gamma$ has entries defined in~\eqref{equ:gammadef}, and thus  $\m \Gamma\bar{\m u}_g$ is an $n \times 1$ vector that gives the generation for each bus of the network. Likewise, matrix $\bm\Pi$ has entries defined in~\eqref{equ:pidef}, and 
$\m \Pi (\m u_{{b}}+\hat{\m w}_{\mathrm{misc}})$ is an $n\times 1 $ vector that gives the building loads per bus. Therefore, vector  $\m \Gamma\bar{\m u}_g -\m \Pi (\m u_{{b}}+\hat{\m w}_{\mathrm{misc}}) - \hat{\m w}_{\mathrm{BL}}$, which appears in both~\eqref{equ:DCOPF2} and~\eqref{equ:DCOPF3}, represents the net nodal power injections. With $\m 1_n \in \mathbb{R}^n $ defined as a vector of all ones, $\left(\m \Gamma\bar{\m u}_g-\m \Pi (\m u_{{b}}+\hat{\m w}_{\mathrm{misc}})- \hat{\m w}_{\mathrm{BL}} \right)^{\top} \m 1_n$ gives the sum of all net nodal injections, and constraint~\eqref{equ:DCOPF2} ensures the supply-demand balance. In addition,   $\m F^{\mathrm{max}} \in \mathbb{R}^{n_l}$ is the vector containing the thermal limits for real power flow on the $n_l$ branches of the network; and $\m L_{\mathrm{ptdf}} \in \mathbb{R}^{n_l \times n}$ is a matrix of power transfer distribution factors \cite{zimmerman2011matpower}. This matrix maps net nodal injections to line power flows, and thus~\eqref{equ:DCOPF3} guarantees the satisfaction of line flow limits.  Formulation~\eqref{equ:DCOPF} is useful in the next sections. }

\section{How Can Buildings Impact Power Grids?\\ Addressing the Main Challenges and BtG-GMPC}~\label{sec:OCP}
In the previous section, we formulate the dynamics of the buildings-integrated power network.  The presence of $\m u_b(t)$  in~\eqref{equ:PlantSS}, exemplifies the control potential that buildings have on power system operation and control, and hence the integration  advocated in this paper. 
In this section, we investigate the discrepancies in time-scales between the building~\eqref{equ:CBuildingDynamincs2} and power network dynamics~\eqref{equ:PlantSS} and discuss a formulation of the joint optimal control problem that addresses the time-scale discrepancies, while seamlessly incorporating objectives and constraints from the power grid and building clusters.  
\subsection{The Not-So-Cruel Curse of Time-Scales}~\label{sec:TS}
The formulated dynamics in Sections~\ref{sec:buildingsSection} and~\ref{sec:PowerNetwork} clearly operate in two different time-scales. While grid regulation problems and mechanical input power variations are often in seconds, the building dynamics and controls are much slower. For example, temperatures in buildings change slowly in comparison with frequencies and voltages in power networks. 

To overcome this limitation, we design local optimal control laws that operate at different scales. Specifically, the time-step for application of building optimal control laws is $h_b$; and the time step for application of grid optimal control laws is $h_g$, where $h_g<<h_b$. This approach reflects the physical realities for these systems, and this consideration can be imposed via constraints in the optimal control problem, whose construction is the objective of this section. 
Since buildings possess slower dynamic behavior, we restrict the controls of buildings to be fixed for the faster time-scale of the power network. 	

Given this integration scheme, the discrepancy in time-scales between building and grid dynamics, the natural existence of algebraic equations in the power network model, and the necessity of including hard constraints such as tight frequency and temperature bounds, model predictive control (MPC) is the natural solution to solve the joint optimal control problem. Other control methods such as rule-based control or PID control can still be used for individual buildings, but these techniques provide inferior results in comparison to MPC as discussed in~\cite{morocsan2011distributed,patel2016distributed,Mirakhorli2016499}. In addition, analytical optimal control techniques that are based on deriving a closed form solution of the optimal control law cannot be computed due to the reasons outlined above. 
\subsection{Discretization of the Dynamics via Gear's Method}
Another challenge with BtG integration is the presence of algebraic equations in~\eqref{equ:PlantSS} emerging from power flows of load nodes. Here, we present a simple, yet high-fidelity discretization routine for two dynamical systems with different time-scales and algebraic constraints.

First, we assume that the sampling times for the power grid [cf.~\eqref{equ:PlantSS}] and building [cf.~\eqref{equ:CBuildingDynamincs2}] dynamics are respectively $h_g$ and $h_b$;
note that $h_b >> h_g$.  The discretization we utilize in this paper is based on Gear's method---a backward differentiation routine---for DAE (descriptor) systems~\cite{Sincovec1981}. The discretization of~\eqref{equ:PlantSS} can be written as follows:\footnote{For the discretization purposes, we use the linearized power flows by assuming that $\sin(\delta_i-\delta_j)=\delta_i - \delta_j$. This implies that the state-space matrix $\m A_g$ in~\eqref{equ:PlantSSd} is different than $\m A_g$ in~\eqref{equ:PlantSS}, where the former includes the Laplacian matrix of the network with weights equal to line inductances. However, we simulate the system with the nonlinear power flows in Section~\ref{sec:results}.}
\begin{eqnarray}
{\m x}_g(k_gh_g) &=& \m f_g(\m x_g,\m u_g,\m u_b, \m w_g)= \bar{\m A}_g \sum_{i=1}^{s} \alpha_{i}\m E_g {\m x}_g(h_g(k_g-i)) \nonumber \\
&&\hspace{-1.85cm}+\m B_0\left({\m A}_{\m u_{b}}  \m u_{b} (k_gh_g) +{\m B}_{u_{g}} \m u_{g}(k_gh_g) +\m B_{w_{g}}\m w_{g}(k_gh_g) \right),~\label{equ:PlantSSd}
\end{eqnarray}
where $\bar{\m A}_g=(\m E_g-h_g\beta_0\m A_g)^{-1}, \m B_0 = h_g\beta_0\bar{\m A}_g, \beta_0 = ({\sum_{i=1}^{s}1/i})^{-1} \text{,} \,$ $\alpha_i = (-1)^{i+1}\beta_0\sum_{j=i}^{s}j^{-1}\binom{j}{i}$; $k_g$ is the time-step for the grid dynamics; {and $s$ is called the order of the method}. This method requires a set of $s$ initial conditions; see Remark~\ref{rem:Gear}.
Similarly, the discrete form of~\eqref{equ:CBuildingDynamincs2} can be written as follows
\begin{eqnarray}
\m x_b(k_bh_b) &&=\m f_b(\m x_b,\m u_b,\m w_b) =\bar{\m A}_b \sum_{i=1}^{s} \alpha_{i} {\m x_b}(h_b(k_b-i)) \nonumber \\
&&+\m B_1\left({\m B}_{u_{b}}  \m u_{b} (k_bh_b) +{\m B}_{w_{b}} \m w_{b}(k_bh_b)  \right),\label{equ:CBBuildingDynamicsDiscrete}
\end{eqnarray}
where $\bar{\m A}_b=(\m I_{2n_b}-h_b\beta_0\m A_b)^{-1}, \m B_1=h_b\beta_0\bar{\m A}_b, \beta_0 = ({\sum_{i=1}^{s}1/i})^{-1} , \alpha_i = (-1)^{i+1}\beta_0\sum_{j=i}^{s}j^{-1}\binom{j}{i}$, and $k_b$ is time-step for dynamic operation of buildings. 

Gear's discretization amounts to a backward Euler-like implicit method. { Gear's method is applicable to DAEs where the matrix $\m E_g$ is allowed to be singular; notice that~\eqref{equ:PlantSSd} does not rely on the inverse of $\m E_g$. Interestingly, when the matrix $\m E_g$ is identity and $s=1$, Gear's method reduces to the standard backward Euler's method.} The principal merit of implicit methods is that they are typically more stable for solving systems with a larger step size $h$, while still performing well for systems with faster time-constants~\cite{Sincovec1981}. 
\begin{rem}[Convergence of Gear's Method]\normalfont \label{rem:Gear}
The states of the discretized descriptor system in~\eqref{equ:PlantSSd} and~\eqref{equ:CBBuildingDynamicsDiscrete} converge to the actual ones in a finite number of time-steps, even if the $s$-initial conditions are arbitrarily chosen~\cite{Sincovec1981}.
A method to compute the \textit{correct} initial conditions is also provided in~\cite{Sincovec1981}. 
\end{rem}
\subsection{Joint Optimal Control Problem: BtG-GMPC}~\label{sec:MPCBuilding}
The joint optimal control problem, Building-to-Grid Gear MPC (BtG-GMPC), is formulated as in~\eqref{equ:JointMPC}. The variables, cost function, and constraints of BtG-GMPC are as follows:
\begin{itemize}[leftmargin=*]
\item $T_p$ is the prediction horizon and $t$ is the initial starting point of the MPC.  The formulation shows the MPC for one prediction horizon.
\item $\m U_b=\{\m u_b(t+h_b),\m u_b(t+2h_b),\ldots,\m u_b(t+T_p)\}, \Delta \m U_g=\{\Delta\m u_g(t+h_g),\Delta\m u_g(t+2h_g),\ldots,\Delta\m u_g(t+T_p)\},$ and $\bar{\m u}_g$ are the three sets of optimization variables that we defined previously. In addition, the two sets of states defined as $\m X_b=\{\m x_b(t+h_b),\m x_b(t+2h_b),\ldots,\m x_b(t+T_p)\}$ and $ \m X_g=\{\m x_g(t+h_g),\m x_g(t+2h_g),\ldots,\m x_g(t+T_p)\}$ are also optimization variables. 
\item The cost function $f(\Delta\m U_g,\bar{\m u}_g,\m U_b, \m X_g, \m X_b)$ is defined as the weighted summation of the building costs, the steady-state LOPF costs, the penalties on the deviation from the steady-state generation, and the deviation cost from the nominal synchronous frequency:
\begin{eqnarray}
	f(\cdot)&=&J(\bar{\m u}_g)+\frac{h_b}{T_p}\sum_{k_b=1}^{{T_p}/{h_b}} \left[ \m c_b^{\top}(t+k_bh_b)\m u_b (t+  k_bh_b)\right]\notag \\
	&& +\dfrac{h_g}{T_p}\sum_{k_g=1}^{{T_p}/{h_g}}  [\Delta \m u_g^{\top} (t+k_gh_g)\m R\Delta \m u_g( t+k_gh_g)  \notag \\
	&& +\m x_g^{\top} (t+k_gh_g)\m Q \m x_g(t+k_gh_g)], \label{equ:BtGobj}
\end{eqnarray}	
where
	\begin{itemize}[leftmargin=*]
\item
	 $J(\bar{\m u}_g)$ is the LOPF cost function~\eqref{equ:DCOPF0}. The parameters of this cost function are widely available in the power systems literature~\cite{zimmerman2011matpower}.
	 \item $\m c_b(t+k_bh_b)$ is a time-varying vector representing the cost of electricity at time $t+k_bh_b$. These prices are the wholesale price of electricity for commercial building operators. 
\item The third term in $f(\cdot)$ penalizes the deviations in the mechanical power setpoints of generators through a quadratic cost function, with matrix $\m R\in \mathbb{R}^{n_g\times n_g}$ being the quadratic penalty matrix, which is assumed to be positive semidefinite.  
\item The fourth term in $f(\cdot)$ penalizes the deviations of the generator frequencies from their nominal value using matrix $\m Q\in \mathbb{R}^{2n \times 2n}$.  The reader is referred to~\cite{Ko2004,Taylor2015} for related constructions.  This cost function is similar to the linear quadratic regulator, which is used in dynamical systems and power network stability studies.  
		\item The terms multiplying the summations are meant to average the building and grid costs across the planning horizon $T_p$. 
	\end{itemize}
	\item Constraints~\eqref{equ:JMPCGridDs}--\eqref{equ:JMPCGridGridBounds2} depict the dynamics of the building-integrated power grid, as well as lower and upper bounds on the states and inputs of the grid states and controls. Note that
	$ {\m x}_g(t+k_gh_g)=\m f_g(\m x_g,\m u_g,\m u_b,\hat{\m w}_g\,|\,t,s)\overset{\Delta}{=} \bar{\m A}_g \sum_{i=1}^{s} \alpha_{i}\m E_g {\m x}_g(t+h_g(k_g-i)) + \m B_0({\m A}_{\m u_{b}}  \m u_{b} (t+k_gh_g) +{\m B}_{u_{g}} \m u_{g}(t+k_gh_g) +\m B_{w_{g}} \hat{\m w}_{g}(t+k_gh_g))$, where $s$ corresponds to the order of Gear's method. 
	\item Constraints~\eqref{equ:JMPCBldgDs}--\eqref{equ:JMPCGridBldgBounds2} represent the building cluster dynamics and the bounds on the states and inputs of the individual buildings, while constraint~\eqref{equ:JMPCLOPF} imposes the constraints of the LOPF as discussed in the previous section.
	\item The final constraint~\eqref{equ:JMPCTimeScales} represents the idea of the time-scales integration whereby the building control variables are kept constant between two consecutive building instances. Since $h_b>h_g$, we assume that between two consecutive building sampling instances (i.e., $k_bh_b$ and $(k_b+1)h_b$), the building controls $u_b(k_bh_b)$ are all constant variables to be found. Hence, for all $\forall \, k_gh_g \in [k_b h_b, (k_b+1)h_b)$, $\m u_b(k_bh_b)=\m u_b(k_gh_g)=\bar{\m u}_b.$ 
\end{itemize}
\begin{mdframed}[style=MyFrame]
	\begin{subequations}~\label{equ:JointMPC}
		\begin{align}
			\hspace{-0.21cm}\textbf{BtG-GMPC:}&\notag\\
			\hspace*{-0.8cm}\minimize_{\substack{\m U_b,\Delta \m U_{g},\bar{\m u}_g\\ \m X_b, \m X_g}}\,\,\, & f(\Delta\m U_g,\bar{\m u}_g,\m U_b, \m X_g, \m X_b) \label{equ:JMPCCost}\\
			\mathrm{subject\,to}\hspace{0.3cm} \; 
			&\hspace{-0.34cm} \,{\m x}_g(t+k_gh_g) = \m f_g(\m x_g,\m u_g,\m u_b,\hat{\m w}_g\,|\,t,s) \label{equ:JMPCGridDs}\\
			& \Delta\m u_g^{\mathrm{min}}\leq  \Delta\m u_g( t+ k_gh_g)\leq \Delta\m u_g^{\mathrm{max}} \label{equ:JMPCGridGridBounds1} \\
			& \m x_g^{\mathrm{min}}\leq  \m x_g(t+k_gh_g)\leq \m x_g^{\mathrm{max}} \label{equ:JMPCGridGridBounds2}\\
			& \forall \; k_g \in \{1,\ldots,T_p/h_g\}\nonumber \\
			&\m x_b(t+k_bh_b) = \m f_b(\m x_b,\m u_b,\hat{\m w}_b\,|\,t,s) \label{equ:JMPCBldgDs} \\
			& \m u_b^{\mathrm{min}}\leq  \m u_b(t+k_bh_b)\leq \m u_b^{\mathrm{max}}\label{equ:JMPCGridBldgBounds1}\\
			& \m x_b^{\mathrm{min}}\leq  \m x_b(t+k_bh_b)\leq \m x_b^{\mathrm{max}} \label{equ:JMPCGridBldgBounds2} \\
			& 	\forall \; k_b \in \{1,\ldots,T_p/h_b\}\nonumber \\
			& \eqref{equ:DCOPF1},\eqref{equ:DCOPF3}\text{\footnotemark}  \label{equ:JMPCLOPF}\\
			& \m u_b(t+k_gh_g)=\bar{\m u}_b=\m u_b(t+k_bh_b)\label{equ:JMPCTimeScales} \\
			& \forall \, k_gh_g \in [k_b h_b, (k_b+1)h_b) \notag . 
		\end{align}
	\end{subequations}
\end{mdframed}
\footnotetext{Constraint~\eqref{equ:DCOPF2} (the supply-demand balance) is removed from BtG-GMPC as it is implicitly present in the discretized algebraic equations in~\eqref{equ:JMPCGridDs}. }
\begin{algorithm}
\small	\caption{Moving Horizon BtG-GMPC \& LOPF Coupling}\label{algo}
	\begin{algorithmic}
\State \textbf{input:} BtG-GMPC forecasts and parameters, $\m x_b(-(s-1)h_b:h_b:0),\m x_g(-(s-1)h_g:h_g:0),T_p,T_{\mathrm{final}}$
\State \textbf{output:} $\{\bar{\m u}^*_g,\Delta \m u_g^*, \m u_b^*\} \;\;\forall t \in [0,T_{\mathrm{final}}]$
\State \textbf{while} $t<T_{\mathrm{final}}$
\State \hspace{0.5cm} \textbf{if} $t=\kappa T_p$ (multiple of $T_p$, i.e., $t=0,T_p,2T_p,\ldots$)
\State \hspace{0.8cm} \textbf{solve}  BtG-GMPC~\eqref{equ:JointMPC} for $\m U_b^*,\Delta \m U_{g}^*,\bar{\m u}_g^*$ 
\State \hspace{0.8cm} \textbf{apply} $\bar{\m u}_g^*$ $\;\forall t \in [\kappa T_p,(\kappa+1)T_p]$
\State \hspace{0.8cm} \textbf{apply} $\m U_b^*(1)$ $\;\forall t \in [t,t+h_b]$
\State \hspace{0.8cm} \textbf{apply} $\Delta \m U_{g}^*(1)$ $\;\forall t \in [t,t+h_g]$
\State \hspace{0.8cm} \textbf{discard} $\m U_b^*(2:\mathrm{end}),\Delta {\m U}_g^*(2:\mathrm{end})$
\State \hspace{0.5cm} \textbf{else if} $(t = \kappa_1 h_g) \wedge (t \neq \kappa_2 T_p) \wedge (t \neq \kappa_3 h_b)$ 
\State \hspace{0.8cm} \textbf{solve}~\eqref{equ:JointMPC} \textbf{without} $\bar{\m u}_g,\m U_b$, while eliminating constraints~\eqref{equ:JMPCBldgDs}--\eqref{equ:JMPCTimeScales}
where $\m U_b,\bar{\m u}_g$ are the optimal constant values from the previous/subsequent steps
\State \hspace{0.8cm} \textbf{apply} $\Delta \m U_{g}^*(1)$ $\;\forall t \in [t,t+h_g]$
\State \hspace{0.8cm} \textbf{discard} $\Delta \m U_g^*(2:\mathrm{end})$
\State \hspace{0.5cm} \textbf{else if} $(t = \kappa_1 h_b ) \wedge (t \neq \kappa_2 T_p)$
\State \hspace{0.8cm} \textbf{solve}~\eqref{equ:JointMPC}  \textbf{without} $\bar{\m u}_g$, \eqref{equ:JMPCLOPF}, and $J(\bar{\m u}_g)$ 
\State \hspace{0.8cm} \textbf{apply} $\m U_b^*(1)$ $\;\forall t \in [t,t+h_b]$
\State \hspace{0.8cm} \textbf{apply} $\Delta \m U_{g}^*(1)$ $\;\forall t \in [t,t+h_g]$
\State \hspace{0.8cm} \textbf{discard} $\m U_b^*(2:\mathrm{end}),\Delta \m U_g^*(2:\mathrm{end})$
\State \hspace{0.5cm}\textbf{end if}
\State \hspace{0.5cm}$t\leftarrow t+h_g$
\State \textbf{end while}
	\end{algorithmic}
\end{algorithm}
Algorithm~\ref{algo} illustrates a routine that implements BtG-GMPC's rolling horizon window along with the integration of the LOPF problem. Given the BtG-GMPC parameters (including the first $s$-initial steps of the discretized dynamics), the algorithm computes the optimal solutions to the LOPF problem and the joint MPC. For simplicity, we assume that the prediction horizon $T_p$ is equivalent to the time-scale in which the optimal dispatch is solved, i.e., 5 to 15 minutes. 

We also assume that $h_g < h_b < T_p << T_{\mathrm{final}}$ and $h_b/h_g$, $T_p/h_b$, $T_p/h_g$ are all positive integers. The algorithm starts by finding the solution to the generator's operating points $\bar{\m u}_g$ for any multiple of the prediction horizon $T_p$, as well as the deviation from this setpoint $\Delta \m u_g(t)$ and $\m u_b(t)$ up until the next planning horizon, and so on. As in classical MPC routines, only the first instance of the optimal control trajectory is applied, while the rest are discarded. Note that the BtG-GMPC with LOPF is only solved for when $t$ (the counter) is a multiple of $T_p$. 
If $t$ is not a multiple of $T_p$, but a multiple of the building's sampling time $h_b$, the building and grid controls are computed. The final case captures the gap between the two time-scales: where the building and grid controls are applied, the building controls are kept constant from the previous optimal computations, while grid controls are computed in the meantime for every grid sampling time.
\begin{rem}[Tractability of BtG-GMPC]\normalfont
Problem~\eqref{equ:JointMPC} is a quadratic program. Even for large-scale systems, this optimization routine is tractable, and can be solved by off-the-shelf solvers such as \texttt{CPLEX}, \texttt{MOSEK}, or Matlab's \texttt{QuadProg}. 
\end{rem}

\begin{rem}[Fast MPC and Time-Complexity]\label{rem:FastMPC}\normalfont
The BtG-GMPC optimization is applied online as predictions for the uncontrollable inputs might not be available for times greater than the prediction horizon $T_p$. However, given that prediction for uncontrollable inputs are available prior to the start of the day, this problem can be solved offline. If solved online, fast online MPC algorithms for quadratic programs have been developed in~\cite{wang2008fast} and can be immediately applied to BtG-GMPC. Otherwise, the problem can be solved offline, which eases the communication requirement of exchanging optimal solutions.
Note that in BtG-GMPC, the maximum total number of variables at each time-step is equal to $3n_b+2n+2n_g=N$. As reported in~\cite{wang2008fast}, MPC formulations take $\mathcal{O}(T_p \cdot N^3)$ at each time-step. This is based on novel interior-point-based implementations.  
\end{rem}

\section{Comparisons with Decoupled BtG Designs}\label{sec:anlytical}

The BtG integration framework developed in this paper enables building and grid operators to jointly optimize their decisions. This section analytically formalizes and compares the BtG framework with decoupled designs, in which the grid and the building operators schedule generation and the building power consumption separately. 

Problem~\eqref{equ:JointMPC} jointly optimizes over two groups of variables: grid decisions $(\bar{\m u}_g, \Delta \m U_g, \m X_g)$ and building decisions $(\m U_b, \m X_b)$. The problem  can be written in the following way, which brings out the coupling between the grid and building decisions.
\begin{subequations}
	\label{equ:JointOptAbstract}
	\begin{align}
	f^*_{\mathrm{BtG}}=\minimize_{\substack{\bar{\m u}_g, \Delta \m U_g, \m X_g, \\ \m U_b, \m X_b}} \;\; & f_g(\bar{\m u}_g, \Delta \m U_g, \m X_g) + f_b(\m U_b, \m X_b) 
	\label{equ:abs-obj}\\
	\subjectto \; \; & (\bar{\m u}_g, \Delta \m U_g, \m X_g, \m U_b) \in \mathcal{C} \label{equ:abs-coupl}\\
	& (\bar{\m u}_g, \Delta \m U_g, \m X_g) \in \mathcal{F}_g \label{equ:abs-gridconstr} \\
	& (\m U_b, \m X_b) \in \mathcal{F}_b. \label{equ:abs-bldgconstr}
	\end{align}
\end{subequations}
The objective~\eqref{equ:abs-obj} corresponds to~\eqref{equ:JMPCCost}. All constraints of problem~\eqref{equ:JointMPC} are captured by one of the constraints of problem~\eqref{equ:JointOptAbstract}. Specifically, set $\mathcal{C}$ represents the coupling between grid and building  decisions, that is, $\mathcal{C}$ is the set of all decisions $(\bar{\m u}_g, \Delta \m U_g, \m X_g, \m U_b, \m X_b)$ that satisfy~\eqref{equ:JMPCGridDs} and~\eqref{equ:DCOPF3}. Set $\mathcal{F}_g$ represents all the constraints pertaining to grid decisions only, that is, constraints~\eqref{equ:DCOPF1}, \eqref{equ:JMPCGridGridBounds1}, and~\eqref{equ:JMPCGridGridBounds2}. Set $\mathcal{F}_b$ represents all constraints pertaining to building decisions only, that is,~\eqref{equ:JMPCBldgDs}, \eqref{equ:JMPCGridBldgBounds1}, \eqref{equ:JMPCGridBldgBounds2}, and~\eqref{equ:JMPCTimeScales}. 

Grid and building controls are optimized separately in traditional power systems. In particular, building operators control the building HVAC load $\m U_b$ based on electricity prices $c_b$, forecasted weather conditions, and occupancy behaviors [with the latter two captured by  $\m w_b(t)$]. Bang-bang control or more sophistacted MPC methods may be used to determine the building HVAC load.  In turn, grid operators  forecast the grid load, which comprises the building load $\m u_b+\hat{\m w}_{\mathrm{misc}}$ and the remaining base load $\hat{\m w}_{\mathrm{BL}}$, and determine generator setpoints and mechanical power adjustments. 

The previously described process can be formalized within the proposed framework as follows. Supposing the building HVAC controls are optimized via MPC, the building operators solve the following optimization problem:
\begin{subequations}
	\label{equ:BldgOptAbstract}
	\begin{align}
	f_b^{\mathrm{MPC}}=\minimize_{\m U_b, \m X_b}\;\;\;& f_b(\m U_b, \m X_b) 
	\label{equ:bldg-abs-obj}\\
\subjectto \;\;\;& (\m U_b, \m X_b) \in \mathcal{F}_b. \label{equ:bldg-abs-bldgconstr}
	\end{align}
\end{subequations}
That is, the objective is to minimize the cost of building operation, subject to the dynamical constraints of the buildings as well as the state and control bounds. It is actually not difficult to see that the previous optimization can be performed by each building operator separately. 
Let $\m U_b^{\mathrm{MPC}}, \m X_b^{\mathrm{MPC}}$ be the solution of problem~\eqref{equ:BldgOptAbstract}.

The grid operator optimizes the generator setpoints and mechanical power adjustments based on the predicted grid load, which includes building loads $\m U_b^{\mathrm{MPC}}$.  The grid operator thus solves the following optimization problem:
\begin{subequations}
	\label{equ:GridOptAbstract}
	\begin{align}
	f_g^{\mathrm{MPC}}=\minimize_{\bar{\m u}_g, \Delta \m U_g, \m X_g}\;\; & f_g(\bar{\m u}_g, \Delta \m U_g, \m X_g)  \label{equ:grid-abs-obj}\\
	\subjectto \;\;& (\bar{\m u}_g, \Delta \m U_g, \m X_g, \m U_b^{\mathrm{MPC}}) \in \mathcal{C} \label{equ:grid-abs-coupl}\\
	& (\bar{\m u}_g, \Delta \m U_g, \m X_g) \in \mathcal{F}_g. \label{equ:grid-abs-gridconstr} 
	\end{align}
\end{subequations}
Let $\bar{\m u}_g^{\mathrm{MPC}}, \Delta \m U_g^{\mathrm{MPC}}, \m X_g^{\mathrm{MPC}}$ be the solution of~\eqref{equ:GridOptAbstract}.

The total cost of operation for the previously mentioned decoupled design is $f_g^{\mathrm{MPC}}+f_b^{\mathrm{MPC}}$, where  $f_g^{\mathrm{MPC}}$ and $f_b^{\mathrm{MPC}}$ are respectively  the optimal values of~\eqref{equ:GridOptAbstract} and~\eqref{equ:BldgOptAbstract}. The relationship with the cost from BtG integration $f^*_{\mathrm{BtG}}$ [cf.~\eqref{equ:JointOptAbstract}] is provided in the following proposition. 
\begin{proposition}
	\label{prop:BtGvcMPC}
	It holds that 
	\begin{equation}
	f^*_{\mathrm{BtG}} \leq f_g^{\mathrm{MPC}}+f_b^{\mathrm{MPC}}.
	\label{equ:BtGvsMPC}
	\end{equation}
\end{proposition}
\begin{IEEEproof}
	Consider the decisions $(\bar{\m u}_g^{\mathrm{MPC}}, \Delta \m U_g^{\mathrm{MPC}}, \m X_g^{\mathrm{MPC}}, \m U_b^{\mathrm{MPC}}, \m X_b^{\mathrm{MPC}})$. These are feasible for problem~\eqref{equ:JointOptAbstract}, because they are feasible for problems~\eqref{equ:GridOptAbstract} and~\eqref{equ:BldgOptAbstract}. Since $f^*_{\mathrm{BtG}}$ is the optimal value of~\eqref{equ:JointOptAbstract}, it holds for any feasible point that
	\begin{equation}
	f^*_{\mathrm{BtG}} \leq  f_g(\bar{\m u}_g^{\mathrm{MPC}}, \Delta \m U_g^{\mathrm{MPC}}, \m X_g^{\mathrm{MPC}}) + f_b(\m U_b^{\mathrm{MPC}}, \m X_b^{\mathrm{MPC}}) 
	\label{equ:BtGvsMPC2}
	\end{equation} 
	But $\m U_b^{\mathrm{MPC}}, \m X_b^{\mathrm{MPC}}$ is the solution of~\eqref{equ:BldgOptAbstract}, and therefore, $ f_b(\m U_b^{\mathrm{MPC}}, \m X_b^{\mathrm{MPC}})=f_b^{\mathrm{MPC}}$ holds. Likewise, $\bar{\m u}_g^{\mathrm{MPC}}, \Delta \m U_g^{\mathrm{MPC}}, \m X_g^{\mathrm{MPC}}$ is the solution of~\eqref{equ:GridOptAbstract}, and $f_g(\bar{\m u}_g^{\mathrm{MPC}}, \Delta \m U_g^{\mathrm{MPC}}, \m X_g^{\mathrm{MPC}}) = f_g^{\mathrm{MPC}}$ holds. Utilizing the latter two optimal values in~\eqref{equ:BtGvsMPC2},~\eqref{equ:BtGvsMPC} follows. 
\end{IEEEproof}

The previous proposition asserts that the decoupled design incurs a total system cost  that is no smaller  than the one of the proposed BtG integration scheme. Intuitively, the BtG design allows to jointly look for grid and building control actions in set $\mathcal{C}$, as opposed to fixing the building controls first, and then solving the  grid optimization problem.

Attention is now turned to the case where bang-bang control is used to determine the HVAC loads. Bang-bang control is the simplest type and most common type of HVAC control where the controller follows a strict temperature set point (e.g., $22.22\degree$C). The HVAC control system is switched on (or off)  as soon as the zone temperature exceeds (or is below) the dead band which is generally $\pm 0.5\degree \mathrm{C}$. Bang-bang control thus does not optimally solve~\eqref{equ:BldgOptAbstract}, because it restricts when the control system is turned on. To make a fair comparison, it is supposed that the dead band is tuned so that the resulting temperatures do not go outside of the intervals specified by~\eqref{equ:JMPCGridBldgBounds2}. Likewise, the control action is not allowed to exceed the bounds dictated by~\eqref{equ:JMPCGridBldgBounds1}. By design, bang-bang control adheres to the building dynamics described by~\eqref{equ:JMPCBldgDs}.

Let $\m U_b^{\mathrm{BB}}, \m X_b^{\mathrm{BB}}$ be the control actions and resulting states of bang-bang control; the previous discussion implies that  $(\m U_b^{\mathrm{BB}}, \m X_b^{\mathrm{BB}})$ is feasible for problem~\eqref{equ:BldgOptAbstract}.  Let $f^{\mathrm{BB}}_b$ be the resulting cost of building operation. Also, let $\check{f}_g^{\mathrm{MPC}}$ be the optimal value of~\eqref{equ:GridOptAbstract} where $\m U_b^{\mathrm{MPC}}$ is replaced by $\m U_b^{\mathrm{BB}}$; the resulting system cost is  $\check{f}_g^{\mathrm{MPC}}+f^{\mathrm{BB}}_b $. The following proposition relates the costs derived from bang-bang control with the costs of MPC based operation.
\begin{proposition}
	\label{prop:BtGvsBB}
	Suppose that $(\m U_b^{\mathrm{BB}}, \m X_b^{\mathrm{BB}})$ is feasible for problem~\eqref{equ:BldgOptAbstract}. Then it holds for the resulting building operation cost that
	\begin{equation}
	f^{\mathrm{MPC}}_b \leq f^{\mathrm{BB}}_b
	\label{equ:bldgMPCvcBB}
	\end{equation}
	and for the system cost that
	\begin{equation}
	f^*_{\mathrm{BtG}}  \leq  \check{f}_g^{\mathrm{MPC}} + f^{\mathrm{BB}}_b.
	\label{equ:BtGvsBB}
	\end{equation}
\end{proposition}
\begin{IEEEproof}
	Eq.~\eqref{equ:bldgMPCvcBB} follows from the fact that $(\m U_b^{\mathrm{BB}}, \m X_b^{\mathrm{BB}})$ is feasible for~\eqref{equ:BldgOptAbstract}, while $f^{\mathrm{MPC}}_b$ is the optimal value of the same problem. 
	
	To prove~\eqref{equ:BtGvsBB}, let $\check{\bar{\m u}}_g^{\mathrm{MPC}}, \Delta \check{\m U}_g^{\mathrm{MPC}}, \check{\m X}_g$ be the solution of~\eqref{equ:GridOptAbstract} where $\m U_b^{\mathrm{MPC}}$ is replaced by $\m U_b^{\mathrm{BB}}$. It follows that   
	$(\check{\bar{\m u}}_g^{\mathrm{MPC}}, \Delta \check{\m U}_g^{\mathrm{MPC}}, \check{\m X}_g^{\mathrm{MPC}}, \m U_b^{\mathrm{BB}}, \m X_b^{\mathrm{BB}})$ is feasible for problem~\eqref{equ:JointOptAbstract}. Thus it holds that
	\begin{align*}
	f^*_{\mathrm{BtG}} \leq  f_g(\check{\bar{\m u}}_g^{\mathrm{MPC}}, \Delta \check{\m U}_g^{\mathrm{MPC}}, \check{\m X}_g^{\mathrm{MPC}}) + f_b(\m U_b^{\mathrm{BB}}, \m X_b^{\mathrm{BB}}) 
	\end{align*}
	from which~\eqref{equ:BtGvsBB} follows.
\end{IEEEproof}
The previous proposition asserts that MPC for building HVAC controls incurs smaller costs than bang-bang control, as long as the bang-bang control adheres to the same constraints as MPC---a fact that has previously been demonstrated in the building literature. But more importantly, similarly to Proposition~\ref{prop:BtGvcMPC}, it is concluded that fixing the building controls to the particular scheme cannot improve the system costs over jointly designing the building and grid controls. 

For easier reference, we refer to the decoupled design where the building HVAC loads are computed via bang-bang control and subsequently the grid is optimized via MPC as Scenario I (with optimal value $\check{f}_g^{\mathrm{MPC}} + f^{\mathrm{BB}}_b$, cf.~Proposition~\ref{prop:BtGvsBB}). The respective design where the building HVAC loads are optimized via MPC is referred to as Scenario II (with optimal value $f_g^{\mathrm{MPC}}+f_b^{\mathrm{MPC}}$, cf.~Proposition~\ref{prop:BtGvcMPC}). The developed BtG framework is referred as Scenario III (with optimal value $f^*_{\mathrm{BtG}}$) . The next section provides numerical simulations that test the previously mentioned designs, and corroborate the analytically derived comparisons.

\section{Case Studies}\label{sec:results}
	In this section, we investigate the impact of the proposed BtG integration on the performance of grid's stability and the cost-effectiveness of building control systems. 
\subsection{Experimental Setup and Parameters}
	The case studies are performed on various power networks for $T_{\mathrm{final}}=24$ hours.  In particular, we use \textit{casefiles} (\texttt{case9,case14,case30,case57}) from Matpower~\cite{zimmerman2011matpower} to test different power networks, parameters, and total number of buildings. Table~\ref{tab:networksetup} documents the total number of buses, generators, buildings  in the aforementioned casefiles.  All of the data and the codes needed to reproduce the results are made available on this page:
	\url{https://github.com/ahmadtaha1/BtG}. The codes are simulated using CPLEX's quadratic program solver~\cite{CPLEX}, and are written in the Matlab environment. The simulations are performed on a PC running Windows 10 Enterprise, Intel(R) Xeon(R) CPU E3-1271 V3 with a 3.60-GHz processor, and 32 GB of RAM. The code allows the testing of any power network with custom-defined total number of buildings at predefined load buses.  
	
		The parameters, exact constraints, weather data, electricity prices, and other details of the problem are carefully chosen to reflect reasonable conditions and IEEE/ASHRAE practices. The parameters are chosen as follows.
		\begin{itemize}[leftmargin=*]
		\item The grid's base load forecast and the miscellaneous loads of all buildings are chosen carefully.  The base-load variations are based on the National Grid New York electricity company's posted demand curve from the Standard Service in New York: \url{https://www9.nationalgridus.com/niagaramohawk/business/rates/5_load_profile.asp}.
		In addition, the disturbances to the buildings, hourly data for $T_{\mathrm{amb}}, \dot{Q}_{\mathrm{sol}}, \dot{Q}_{\mathrm{int}}, P_{\mathrm{BL}},$ and $P_{\mathrm{misc}}$ are all included in the Github link. 
		\item Building loads are modeled based on one reference commercial building located at the main campus of the University of Texas at San Antonio. For this reference building, the construction materials are known and further determined by ASHRAE standard 90.1-2016. The building size is calculated based on design documents.
		\item The RC-parameters for all buildings are obtained using a normal distribution around the following reference (mean) building. The mean parameters are $R_1=R_2=1.16 \times 10^{-4}$ ($\mathrm{\degree C/W}$), $R_{\mathrm{win}}=6.55\times 10^{-3}$ ($\mathrm{\degree C/W}$), $C_{\mathrm{zone}}=7.033 \times 10^{9}$ ($\mathrm{J/\degree C}$), and $C=1.133 \times 10^{9}$ ($\mathrm{J/\degree C}$). The average building size is around $10,000\,\mathrm{m}^2$.
		\item The prices of electricity for the HVAC loads are reproduced from~\cite{patel2016distributed}. The cost functions for generator mechanical power setpoints $\bar{u}_g$ are extracted from Matpower \cite{zimmerman2011matpower}. The same quadratic cost is used for the variations ($\Delta u_g(t)$). The deviation in frequency is penalized with $Q_{k}=50000$ $\$/(\mathrm{rad/sec})^2$ ($Q_k$ is the $k$th diagonal entry in $\m Q$), and the angles are left without any penalties in the $\m Q$-matrix. 
		\item The parameters of the power network, including line parameters $b_{kj},M_k$, and $D_k$ for all buses $k \in \mathcal{B}$, are obtained from Matpower \cite{zimmerman2011matpower} and the power system toolbox~\cite{chow1992toolbox}.
		\item We choose $h_g=10 \sec$, $h_b=300 \sec$, and a prediction horizon $T_p=900 \sec$. For simplicity, 1st order Gear's method is used in the simulations. The value of $h_g$ is consistent with the discrete time interval at which automatic generation control commands are dispatched~\cite[Sec.~12.3]{gloverpowerbook}.
		\item The bound-constraints in~\eqref{equ:JointMPC} are as follows: (a) $ 59 \leq f_k=\omega_k^{\mathrm{true}}/2\pi \leq 61$ (Hz), (b) $21.5 \leq T_{\mathrm{zone}} \leq 23$ (\degree C) for time-periods between 8AM and 8PM, (c) $22 \leq T_{\mathrm{zone}} \leq 25$ (\degree C) for time-periods between 8PM and 8AM, and (d) $0 \leq P_{\mathrm{HVAC}_{l}} \leq 800$ (KW). The limits on the output power of generators can be found in Matpower \cite{zimmerman2011matpower}.
	\end{itemize}
	\begin{table}[t]
		\centering
		\caption{Building-Integrated Power Network Setup.}
		\begin{tabular}{lrrrr}
			\multicolumn{1}{c}{{{}}} & \multicolumn{1}{c}{{{Case 9}}} & \multicolumn{1}{c}{{{Case 14}}} & \multicolumn{1}{c}{{{Case 30}}} & \multicolumn{1}{c}{{{Case 57}}} \\
			\midrule
			Number of Buses & 9     & 14    & 30    & 57 \\
			Number of Generators & 3     & 5     & 6     & 7 \\
			Number of Buildings & 965	&1058	&376	& 1822
			\\
			\bottomrule
			\bottomrule
		\end{tabular}%
		\label{tab:networksetup}%
	\end{table}%
	\subsection{Impact on Frequency Regulation \& Energy Savings }
	In this section, we present the numerical results for the BtG-GMPC and Algorithm~\ref{algo} (also named Scenario III in the sequel), in comparison with solving the optimal control of buildings and power grids separately via MPC (Scenario II). In addition, we compare the results of BtG-GMPC with bang-bang control of HVAC systems, which is still very common in today's industries, combined with grid-only MPC (Scenario I). The three scenarios are analytically discussed in Section~\ref{sec:anlytical}, have increasing sophistication, but use the same parameters, initial conditions, constraints, and costs. The reader is referred to Section~\ref{sec:anlytical} for a comparison between the decoupled designs (Scenarios I and II) and BtG-GMPC. 
	
	For brevity, we only show the plots for \texttt{case57} from~\cite{chow1992toolbox} with 1822 commercial buildings, but present a cost comparison for all other casefiles. The uploaded Github codes contain the data for other simulations with the corresponding figures for building and power network states and optimal control inputs.  In addition, Table~\ref{tab:comp} shows the cost comparison between the different scenarios. The cost functions are defined as follows: $\m x_g^{\top}\m Q \m x_g$ denotes the frequency deviation cost;  $\Delta \m u_g^{\top} \m R\Delta \m u_g$ represents the mechanical input power deviations cost; $J(\bar{\m u}_g)$ is the LOPF cost; $ \m c_b^{\top} \m u_b$ depicts the HVAC cooling load costs. Note that these costs are all multiplied by \$1,000. In the next section, we compute the perfect cost reduction and compare the different scenarios. The percent reduction in cost is computed as follows:
	\begin{equation}
	\% \; \text{reduction} = \frac{\text{(Cost in Scenario X)} - \text{(Cost in Scenario Y)}}{\text{(Cost in Scenario X)}}
	\notag
	\end{equation}
	where the `Cost' refers to any of the reported cost functions in Table~\ref{tab:comp} and X,Y correspond to any of the three scenarios. 
	\subsubsection{Scenario I} 
\begin{figure*}
	\centering  
	\subfigure[BtG-GMPC: Bus fequencies and optimal power generation.]{\includegraphics[width=0.475\linewidth]{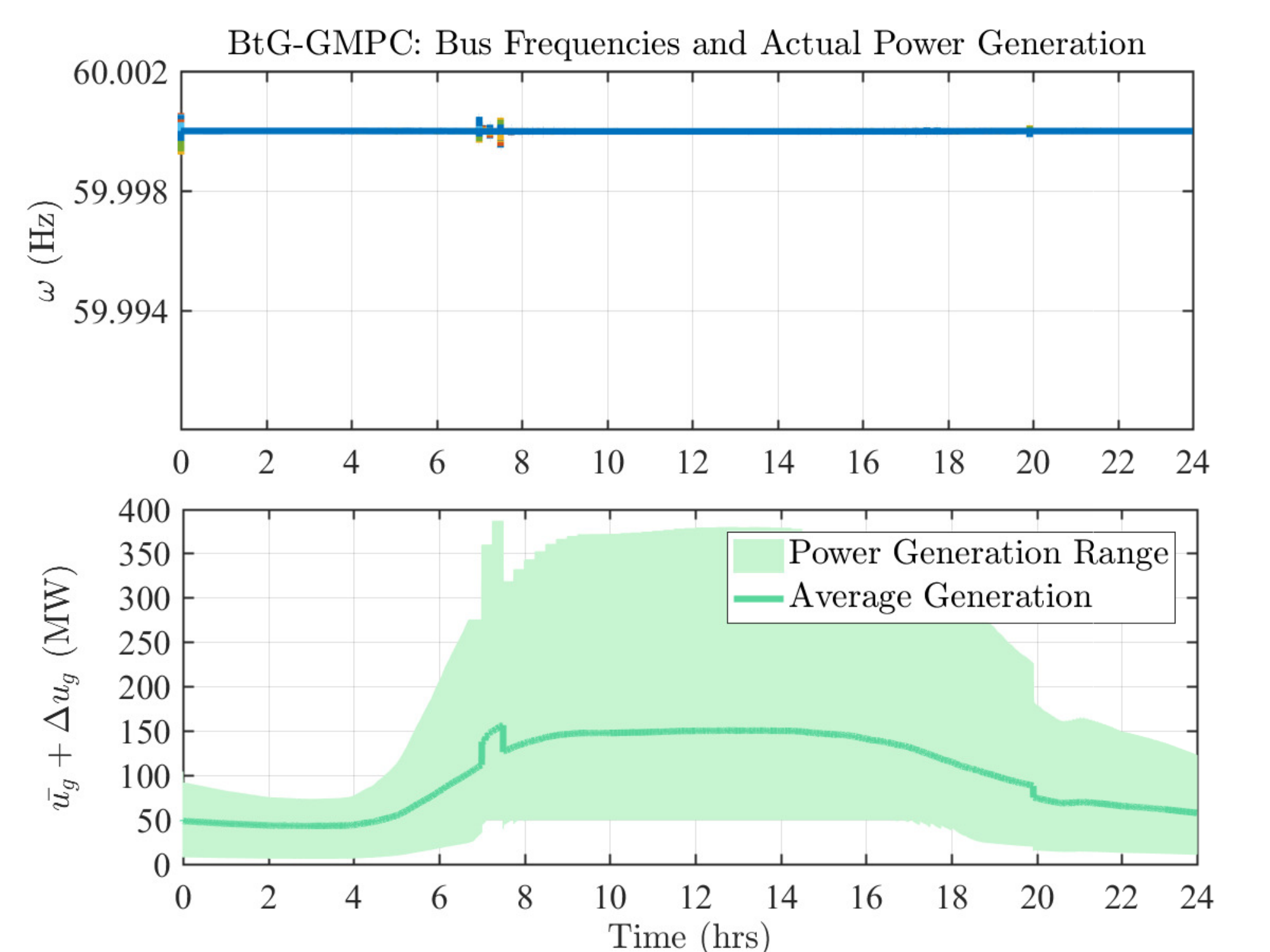}}
	\subfigure[BtG-GMPC: HVAC power consumption and zone temperatures.]{\includegraphics[width=0.475\linewidth]{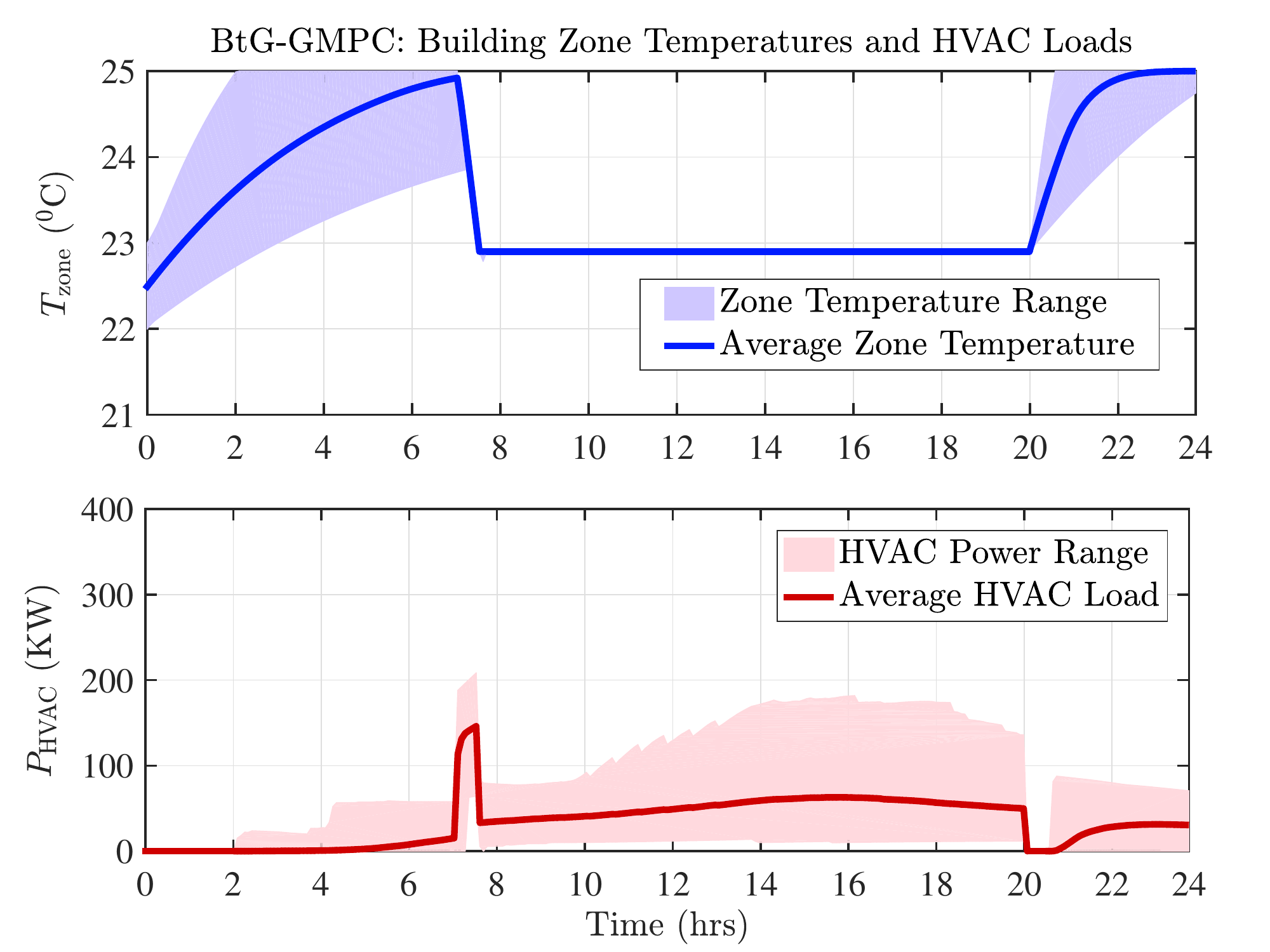}}
	\subfigure[Scenario I: Bus fequencies and optimal power generation.]{\includegraphics[width=0.475\linewidth]{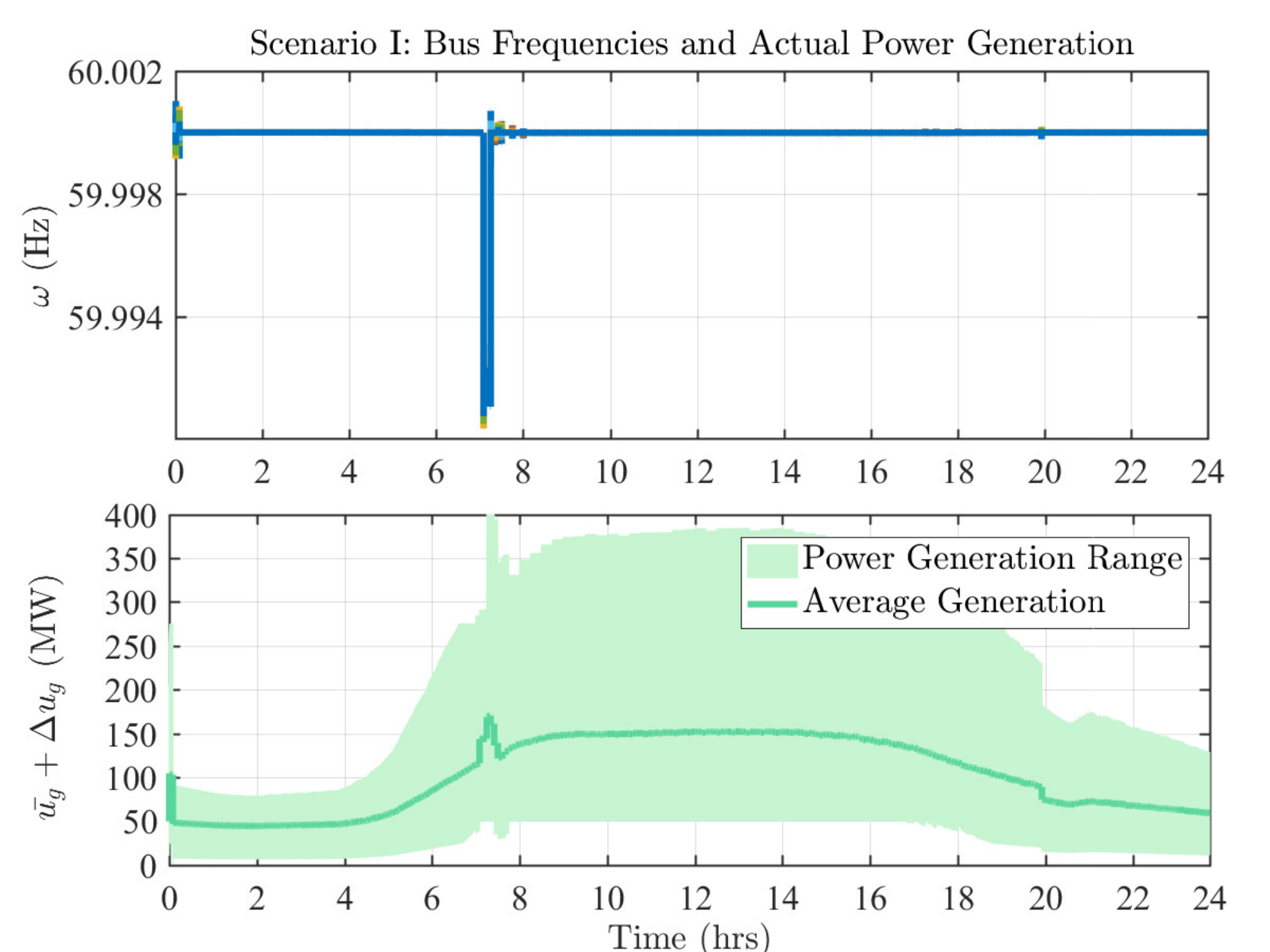}}
	\subfigure[Scenario I: HVAC power consumption and zone temperatures.]{\includegraphics[width=0.475\linewidth]{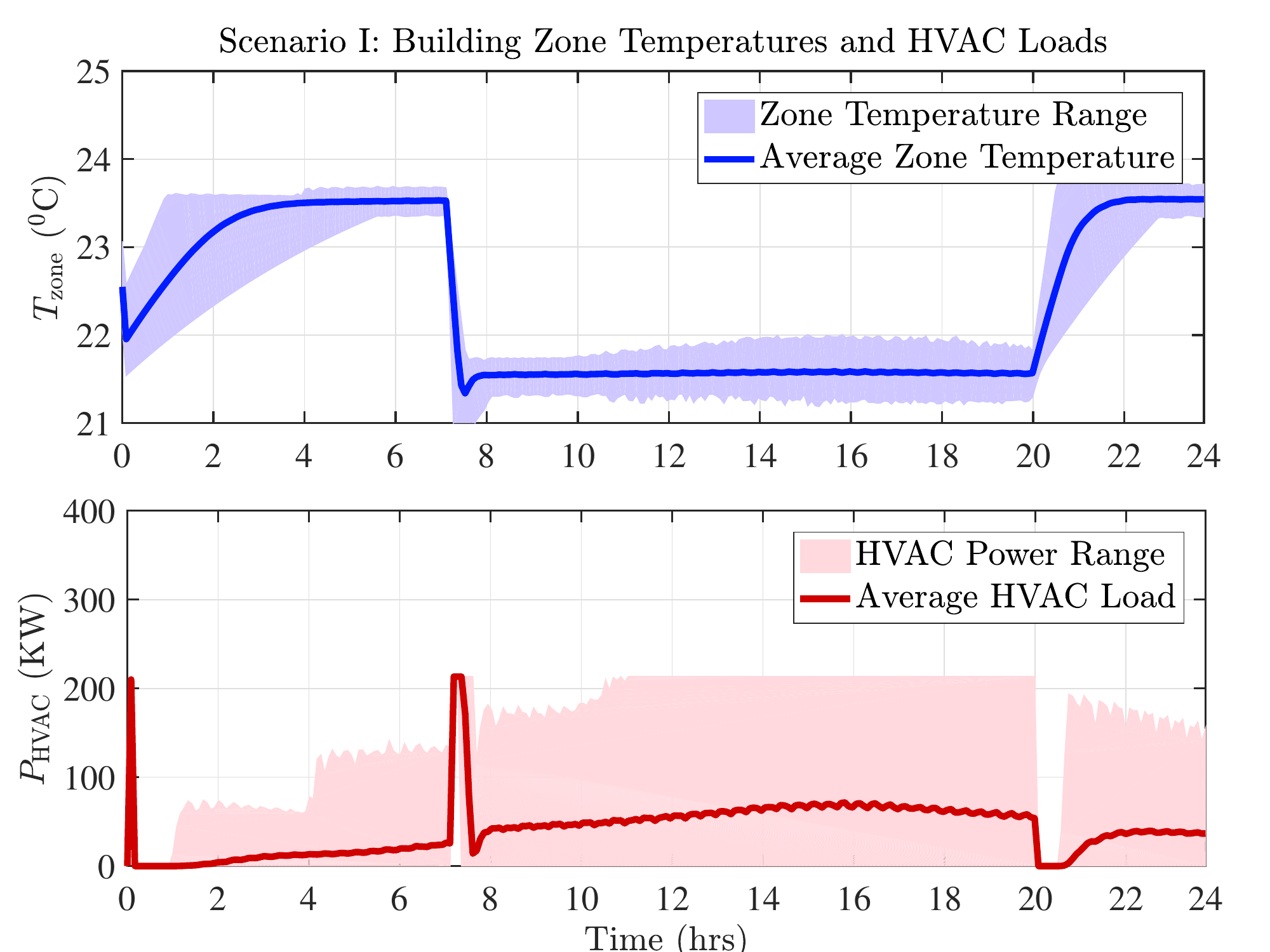}}
	\caption{Numerical results for Scenarios I and III simulated on \texttt{case57} with 1822 buildings. The figures show the power generation range and average for all generators, frequencies ($\omega$ or $f$ in Hz) of all the buses, HVAC power range and average for all buildings, and the range and average of zone temperatures. }
	\label{fig:simulations}
\end{figure*}

The resulting HVAC power consumptions and corresponding zone temperatures for all buildings are shown in Fig.~\ref{fig:simulations}-(c,d). As expected, the bang-bang building control maintains the temperature in the aforementioned band.

	After simulating this case over a period of 24 hours, these bang-bang HVAC inputs are provided to the grid optimal control problem, which is formulated as an MPC based on Gear's method. Figs.~\ref{fig:simulations}-(c,d) show the frequencies, zone temperatures, total generation and HVAC loads for all 1822 buildings and 57 buses in the network. Due to the intermittent nature of the HVAC load, the grid frequency experiences significant deviations from its nominal value (60 Hz); see Fig.~\ref{fig:simulations}-(c). The frequency variations are more prominent at durations when the total load is at peak values. The total costs for Scenario I are provided in Table~\ref{tab:comp}.
	\begin{table}[t]
	\centering \caption{Cost Comparison in 1,000\$ for the Decoupled Problems (Scenarios I and II) versus BtG-GMPC. }
		\renewcommand{\arraystretch}{1.3}
	    \begin{tabular}{l|llll}
		Test 
		Case & \multicolumn{1}{c}{Cost Function} & \multicolumn{1}{c}{Scenario I} & \multicolumn{1}{c}{Scenario II} & \multicolumn{1}{c}{BtG-GMPC} \\
		\midrule
		\multirow{6}[1]{*}{Case 9} & $\sum \m x_g^{\top}\m Q \m x_g$ & 490.51 & 205.57 & 20.33 \\
		&  $\sum \Delta \m u_g^{\top} \m R\Delta \m u_g$ & 60.00 & 58.95 & 59.39 \\
		& $\sum J(\bar{\m u}_g)$ & 125.08 & 121.05 & 120.26 \\
		& $\sum \m c_b^{\top} \m u_b $ & 656.73 & 548.40 & 548.56 \\
		& Total Grid Cost & 675.59 & 385.56 & 199.98 \\
		& \textit{\textbf{Total Cost}} & \textit{\textbf{1332.33}} & \textit{\textbf{933.96}} & \textit{\textbf{748.54}} \\
				\midrule
		\multirow{6}[0]{*}{Case 14} & $\sum \m x_g^{\top}\m Q \m x_g$ & 289.19 & 134.75 & 3.02 \\
		&  $\sum \Delta \m u_g^{\top} \m R\Delta \m u_g$ & 70.49 & 69.60 & 70.39 \\
		& $\sum J(\bar{\m u}_g)$ & 81.32 & 77.61 & 77.22 \\
		& $\sum \m c_b^{\top} \m u_b $ & 705.48 & 588.12 & 588.24 \\
		& Total Grid Cost & 441.01 & 281.96 & 150.62 \\
		& \textit{\textbf{Total Cost}} & \textit{\textbf{1146.49}} & \textit{\textbf{870.08}} & \textit{\textbf{738.86}} \\
				\midrule
		\multirow{6}[0]{*}{Case 30} & $\sum \m x_g^{\top}\m Q \m x_g$ & 49.49 & 20.69 & 0.02 \\
		&  $\sum \Delta \m u_g^{\top} \m R\Delta \m u_g$ & 3.19  & 3.15  & 3.12 \\
		& $\sum J(\bar{\m u}_g)$ & 1.44  & 1.31  & 1.33 \\
		& $\sum \m c_b^{\top} \m u_b $ & 253.02 & 211.30 & 211.31 \\
		& Total Grid Cost & 54.12 & 25.14 & 4.46 \\
		& \textit{\textbf{Total Cost}} & \textit{\textbf{307.14}} & \textit{\textbf{236.45}} & \textit{\textbf{215.78}} \\
			\midrule
		\multirow{6}[1]{*}{Case 57} & $\sum \m x_g^{\top}\m Q \m x_g$ & 466.91 & 25.75 & 0.09 \\
		&  $\sum \Delta \m u_g^{\top} \m R\Delta \m u_g$ & 159.61 & 156.95 & 155.40 \\
		& $\sum J(\bar{\m u}_g)$ & 49.31 & 46.40 & 46.48 \\
		& $\sum \m c_b^{\top} \m u_b $ & 1230.20 & 1027.18 & 1027.18 \\
		& Total Grid Cost & 675.83 & 229.10 & 201.97 \\
		& \textit{\textbf{Total Cost}} & \textit{\textbf{1906.03}} & \textit{\textbf{1256.28}} & \textit{\textbf{1229.15}} \\
		\bottomrule
		\bottomrule
	\end{tabular}
		\label{tab:comp}
	\end{table}%
\subsubsection{Scenario II} In this scenario, we solve the MPC problems for buildings and the power grid separately. First, the building optimal controls are computed via the same MPC formulated in~\eqref{equ:JointMPC}, while eliminating the power grid constraints and variables. The MPC solution for the building's HVAC loads is then fed into a grid-only MPC. This scenario is useful in the sense that grid operators can model the building's load via a classical building MPC model---this can be viewed as a \textit{decoupled} BtG-GMPC. The numerical results of this case show an improvement in the HVAC power consumption of buildings from Scenario I. Specifically, the cooling load costs decreased (from Scenario I to Scenario II) by 42.9, 16.6, 16.4, and 16.5\% for Case 9, Case 14, Case 30, and Case 57. In addition, the overall system costs decreased by an average of 29.7\% for all casefiles (from Scenario I to II). 
	\subsubsection{Scenario III (BtG-GMPC)} In this scenario, we test the performance of Algorithm~\ref{algo} and the corresponding optimization problem~\eqref{equ:JointMPC}. Scenario I and II are separately compared to Scenario III.
	\begin{table}[t]
		\centering
		\caption{Computational Time for Running BtG-GMPC and Algorithm~\ref{algo}.}
		\begin{tabular}{lllll}
			\multicolumn{1}{c}{\textit{\textbf{ }}} & \multicolumn{1}{c}{{{Case 9}}} & \multicolumn{1}{c}{{{Case 14}}} & \multicolumn{1}{c}{{{Case 30}}} & \multicolumn{1}{c}{{{Case 57}}} \\
			\midrule
			Computational Time (Hours)  &0.26&	0.50&	0.82&	3.74
			\\
			\bottomrule
			\bottomrule
		\end{tabular}%
		\label{tab:computationaltime}
	\end{table}
As asserted by the analytical discussion in Section~\ref{sec:anlytical}, the BtG-GMPC results show significant improvement in grid's frequency deviations and the overall costs as shown in Fig.~\ref{fig:simulations}-(a) and Table~\ref{tab:comp}. The frequency variations are notably lower in comparison with the previous two scenarios, and this is clear from a comparison between Figs.~\ref{fig:simulations}-(a) and (c). Also, there is a significant reduction in the overall cost of operation. The results show around 43, 35, 30, and 36\% total cost reduction between Scenarios I and III (for Case 9, Case 14, Case 30, and Case 57), and 20, 15, 9, and 2\% total cost reduction between Scenarios II and III. Also, the grid frequency deviations are reduced,  leading to a decrease in the cost of grid operation by an average of 74.53\% from Scenario I to III and 46.75\% from Scenario II to III for all casefiles. Finally, the plots in Fig.~\ref{fig:simulations} show that the zone temperatures are well-maintained within the required range, with very little fluctuation (in comparison with Scenario I), and the total HVAC consumption is reduced by an average of 16.52\% which is also clear from Figs.~\ref{fig:simulations}-(b,d). 
For brevity, we do not show the plots for Scenario II and other scenarios for other power networks as they are included in the Github link. 

\subsubsection{Computational Speed}
The computational times for running a 24-hour time-horizon are shown in Table~\ref{tab:computationaltime} for the four power networks. The power network \texttt{Case57} requires 3.74 hours to run for the entire 24-hour simulation horizon. The computational times project that if a larger power network or a micro-grid is considered, with potentially 10,000 buildings or more (and thousands of buses), Algorithm~\ref{algo} can still be implemented in real time. This is due to two reasons. Firstly, a system operator or a utility company (see Remark~\ref{rem:WhoSolves?}) will have more computational power at their disposal. Secondly, a fast MPC routine as the one described in Remark~\ref{rem:FastMPC} can also be implemented to speed up the computations.

\subsubsection{Robustness to Forecast and Model Uncertainty: Offline, Day-Ahead Solutions} 

In the previous section, Fig.~\ref{fig:simulations} and Table~\ref{tab:comp} present the results from the MPC solution ($\m X_g, \m X_b, \m U_b, \m U_g$) using the predicted unknown inputs ($\hat{\m w}_g, \hat{\m w}_b$) for the linearized power system using Gear's method. In this section, we 1) extract the optimal MPC control variables $(\m U_b, \m U_g)$ for the entire simulation horizon; 2) feed these inputs to the nonlinear continuous-time DAE solver for~\eqref{equ:PlantSS} and linear ODE solver for~\eqref{equ:BuildingSS}; 3) add zero-mean Gaussian noise with 10\% standard deviation from the unknown inputs $\hat{\m w}_g$ in~\eqref{equ:PlantSS} and $\hat{\m w}_b$ in~\eqref{equ:BuildingSS}; and 4) perturb the 3R-2C building model~\eqref{equ:BuildingSS} with zero-mean Gaussian noise with 10\% standard deviation from the nominal matrices $\m A_b, \m B_{u_{b}}$, and $\m B_{w_{b}}$. 
	
Four major reasons justify this numerical simulation, namely, 1) to validate Gear's discretization method; 2) to assess the performance of the integration framework on the nonlinear continuous-time DAE model for the grid under mismatch between the forecasted and true disturbances; 3) to test the temperature behavior of buildings under model uncertainties that could be a result of parametric misidentification from building operators; and 4) to examine a scenario where the demand response signals are communicated a day prior to the schedules. The latter essentially alleviates the real-time communication burden of the BtG integration framework, by allowing the grid operators to send the demand response schedules way ahead in time and thereby avoiding the necessity to communicate the schedules in real time.

\begin{figure}[t]
	\centering 	\includegraphics[scale=0.3]{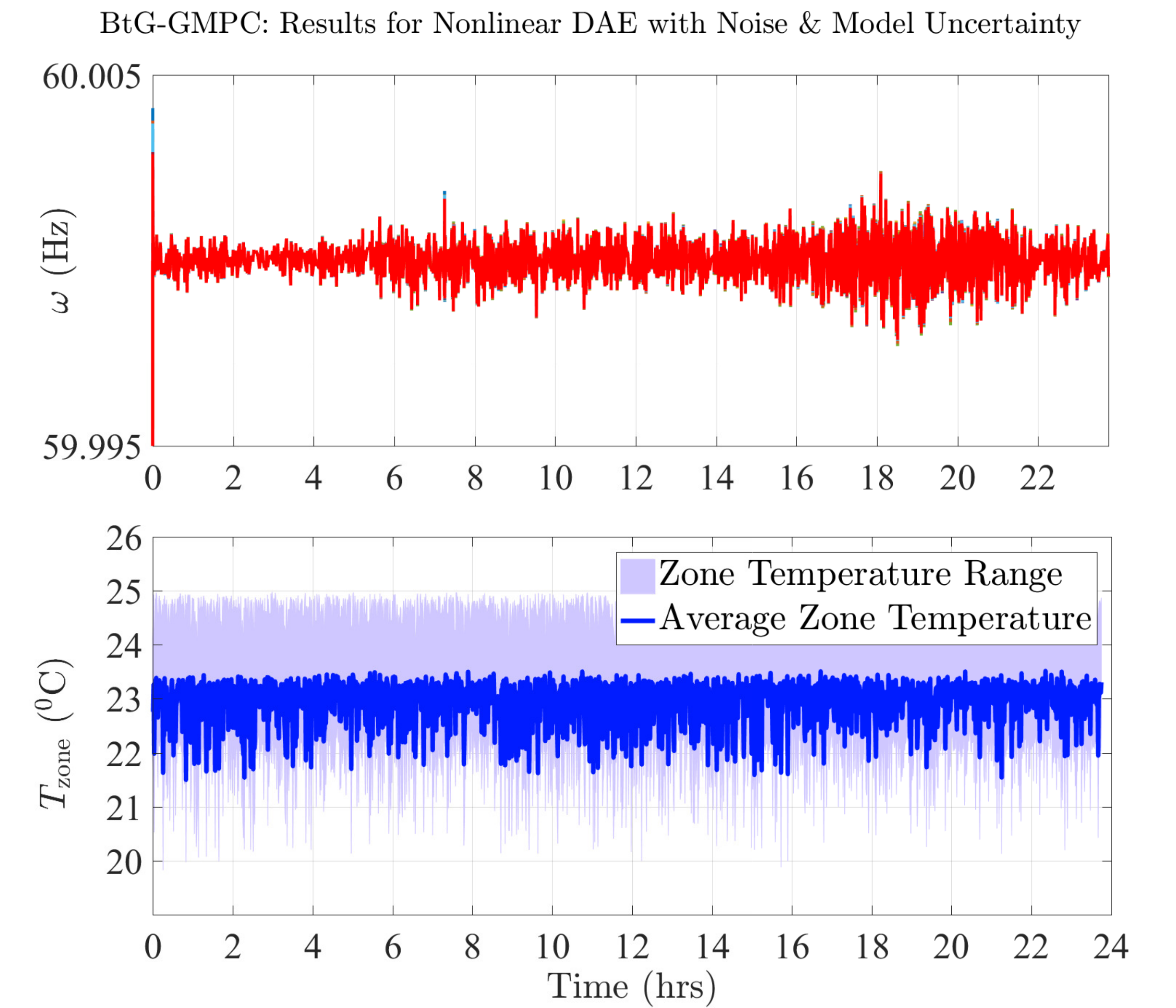}
		\caption{Performance of BtG-GMPC under Forecast and Model Uncertainty: This simulation is performed using a nonlinear DAE solver  for the grid dynamics (Matlab's \texttt{ode15i}) and ODE solver for the building dynamics. The plots depict the grid frequency and zone temperatures for all buildings.}
		\label{fig:NLDAENoise}
	\end{figure}

Fig.~\ref{fig:NLDAENoise} illustrates that the frequency deviations and zone temperatures are kept within reasonable ranges, even under significant parametric and load uncertainty (10\% load mismatch is relatively large in power networks), and the results hence depict that BtG-GMPC is robust to significant disturbances. Note that the MPC scheme is solved offline given the prediction of the loads and the temperatures, and the resulting controls were used as inputs to the nonlinear DAE solver, which demonstrated the good performance of the developed framework. Specifically, the zone temperatures for all buildings are still within the acceptable range, although some buildings experience zone temperatures of $25^{\circ}$C and $20^{\circ}$C, which is due to the parametric mismatch. This can be compared to Fig.~\ref{fig:simulations}-b. In addition, the load mismatch does not destabilize the grid's frequencies; see Fig.~\ref{fig:NLDAENoise} and Fig.~\ref{fig:simulations}-a for comparison. 
	\begin{rem}[Who Solves BtG-GMPC?]~\label{rem:WhoSolves?}\normalfont
BtG-GMPC assumes the knowledge of various parameters such as building RC-constants and generator cost curves. We consider the following: 1) A system operator or a large utility ideally solves BtG-GMPC; 2) commercial building operators contributing to this routine are required to provide modeling parameters for their buildings; and 3) the global signals computed are communicated to the now-contributing operators of individual buildings and generators. The added value of this coupling is two-fold. First, the theoretical impact of expanding the feasible space of two separate problems ensures that the coupled problem's solution will be superior to the decoupled one as illustrated in the previous section. Second, the coupling translates into tangible impact for buildings and the grid, as has been demonstrated in this section. With that in mind, the communication costs are not considered here, and it is assumed that the computed optimal setpoints are communicated instantly to individual buildings and generators.	
\end{rem}
 
\section{Paper Summary and Limitations; Future Work}\label{sec:summary}

This paper introduces the first explicit building-to-grid integration dynamical model with optimal power management formulations and different time-scales.
The paper considers realistic, high-level building models and frequency-focused grid dynamics, in addition to algebraic equations modeling the nodes without generation. We also introduced Gear's method as a high fidelity DAE discretization routine that is leveraged to model BtG integration. The developed framework and optimization problem BtG-GMPC provides setpoints for individual buildings and power grid generators, as well as buildings-aware, optimized optimal power flow setpoints for generators. The formulated problem can be solved efficiently using any quadratic program solver.
Case studies have demonstrated the impact of BtG-GMPC on reducing overall energy costs and minimizing frequency deviations. 

We have kept the dynamical models simple as the focus is on energy consumption and frequency deviations. However, the framework is general and interested researchers can seamlessly extend BtG-GMPC to include advanced models.  Given that, three main challenges are not addressed here. 
\begin{itemize}[leftmargin=*]
	\item We do not consider the problem of controlling or adjusting the reactive power of buildings, and its impacts on regulating the grid's voltages. 
Since the focus of this work is on the framework with specific impacts on frequency regulation, we leave this natural extension to future work. Note that Gear's method and the discretization still holds for models with reactive powers and voltages. 
	\item   The BtG framework developed in this paper pertains to normal grid operation. Unplanned incidents and contingencies  can occur and require appropriate response. In order to deal with plausible contingencies, system operators include reserve scheduling in OPF~\cite{Galiana-ProcIEEE}. Incorporating reserve scheduling is thus an interesting future direction.
 	\item The impact of slow communications between grid and building operators,  load prediction errors, and mismatch in building parameters are all investigated in Section~\ref{sec:results}, and shown to have little impact on the system states.  A more sophisticated BtG-GMPC~\eqref{equ:JointMPC} that incorporates uncertainty in loads and building models is an important improvement that yields a BtG routine tolerant to these unknown inputs.
\end{itemize}
\bibliographystyle{ieeetran}
\bibliography{bibfile}
\begin{IEEEbiographynophoto}{Ahmad F. Taha}
	(S'07--M'15) received the B.E. and Ph.D. degrees in Electrical and Computer Engineering from the American University of Beirut, Lebanon in 2011 and Purdue University, West Lafayette, Indiana in 2015. In Summer 2010, Summer 2014, and Spring 2015 he was a visiting scholar at MIT, University of Toronto, and Argonne National Laboratory. Currently he is an assistant professor with the Department of Electrical and Computer Engineering at The University of Texas, San Antonio since 2015. Dr. Taha is interested in understanding how complex cyber-physical systems operate, behave, and \textit{misbehave}. His research focus includes optimization and control of power system, buildings-to-grid integration, observer design and dynamic state estimation, cyber-security of dynamic systems, and sensor and actuator selection methods. 
\end{IEEEbiographynophoto}

\begin{IEEEbiographynophoto}{Nikolaos Gatsis} received the Diploma degree in Electrical and Computer Engineering from the University of Patras, Greece, in 2005 with honors. He completed his graduate studies at the University of Minnesota, where he received the M.Sc. degree in Electrical Engineering in 2010, and the Ph.D. degree in Electrical Engineering with minor in Mathematics in 2012. He is currently an Assistant Professor with the Department of Electrical and Computer Engineering at the University of Texas at San Antonio. His research interests lie in the areas of smart power grids, renewable energy management, communication networks, and cyber-physical systems, with an emphasis on optimal resource management. Prof. Gatsis co-organized symposia in the area of Smart Grids in IEEE GlobalSIP 2015 and IEEE GlobalSIP 2016. He also served as a Technical Program Committee member for symposia in IEEE SmartGridComm 2013 through 2017.
\end{IEEEbiographynophoto}

\begin{IEEEbiographynophoto}{Bing Dong} received the M.S degree in Building Science from National University of Singapore (NUS), and Ph.D. Degree from Carnegie Mellon University in 2010.  He worked at United Technologies Research Center as Senior Research Scientist from 2010 to 2012. In January 2013, he joined the faculty of Mechanical Engineering at the University of Texas at San Antonio as assistant professor. His research interest is intelligent building operations, modeling and simulation of occupancy behavior in buildings, development of control strategies for buildings-to-grid integration and apply machine learning technologies for future smart building and communities. His research is supported by National Science Foundation, Department of Energy, and industries. 
	
\end{IEEEbiographynophoto}

\begin{IEEEbiographynophoto}{Ankur Pipri} received his B.E. in Electrical Engineering from the National Institute of Technology Bhopal, India and a Bachelor of Technology in October 2012. He then pursued his career as an Operation and Efficiency Engineer in a power plant in Amravati, India. After acquiring the industrial experience, he pursued a Master's degree in Electrical and Computer Engineering from The University of Texas at San Antonio, and graduated in August 2017. His research interest areas include smart grid operations, optimal control methods, and distributed energy resources.
\end{IEEEbiographynophoto}

\begin{IEEEbiographynophoto}{Zhaoxuan Li} received his B.S. in Mechanical Engineering from Shanghai Jiaotong University. (2010), and a M.E. from Texas Tech University (2013). He is currently a Ph.D. candidate in Mechanical Engineering at the University of Texas at San Antonio. He has been involved in urban scale modeling and simulation of smart buildings and smart grids. 
	
\end{IEEEbiographynophoto}

\end{document}

%% file: preamble.tex
 \usepackage{multirow}
 \usepackage{booktabs}
\usepackage{graphics} 

\usepackage{algorithm}
\PassOptionsToPackage{noend}{algpseudocode}
\usepackage{algpseudocode}

\usepackage{gensymb}
\usepackage[colorlinks = true,linkcolor = blue,urlcolor  = blue,citecolor = blue,anchorcolor = blue]{hyperref}

\usepackage{amssymb}

\usepackage{amsthm}
\usepackage{amsmath} 
\usepackage{amssymb}  
\allowdisplaybreaks[4]
\usepackage{subfigure}
\newcommand{\m}{\mathbf}
\usepackage{tikz}
\usepackage{lipsum} 
\usepackage[framemethod=TikZ]{mdframed}
\mdfdefinestyle{MyFrame}{%
	linecolor=black,
	outerlinewidth=1pt,
	roundcorner=4pt,
	innerrightmargin=5pt,
	innerleftmargin=5pt,}

\usepackage{pgfplots}

\usetikzlibrary{decorations.pathmorphing} 
\tikzstyle{block} = [draw,rectangle,thick,minimum height=3em,minimum width=5em]
\tikzstyle{branch} = [circle,inner sep=1pt,minimum size=2mm,fill=white,draw=black]
\tikzstyle{branch1} = [inner sep=1pt,minimum size=0.1mm,fill=white,draw=black]
\tikzstyle{branch1} = [circle,inner sep=0pt,minimum size=1mm,fill=black,draw=black]

\tikzstyle{connector2} = [-,thick]
\tikzstyle{connector} = [->,thick]
\tikzstyle{connector1} = [->,ultra thick]
\tikzstyle{connector4} = [-,ultra thick]

\usepackage[noadjust]{cite}
\newtheorem{theorem}{Theorem}

\newtheorem{rem}{Remark}

\newtheorem{proposition}[theorem]{Proposition}

\DeclareMathOperator{\diag}{diag}

\DeclareMathOperator*{\minimize}{minimize}

\DeclareMathOperator{\subjectto}{subject\ to}

\newcommand{\itn}[1]{^{(#1)}}

\newcommand{\bmat}[1]{\begin{bmatrix} #1 \end{bmatrix}}

\newcommand{\bm}{\boldsymbol}

\usepackage{enumitem}